\documentclass{amsart}

\usepackage{latexsym}
\usepackage{amssymb}
\usepackage{amsmath}
\usepackage{graphicx}

\newtheorem{THM}{Theorem}
\newtheorem{Lemma}[THM]{Lemma}

\newtheorem{cor}[THM]{Corrolary}

\renewcommand{\Re}{\mathbb{R}}

\newcommand{\be}{\begin{equation}}
\newcommand{\ee}{\end{equation}}

\newcommand{\mS}{\mathbb{S}}
\newcommand{\Zed}{\mathbb{Z}}

\newcommand{\Rn}{\Re^{n}}

\newcommand{\ol}[1]{\overline{#1}}
\newcommand{\ra}{\rightarrow}
\newcommand{\da}{\searrow}

\newcommand{\ip}[2]{\left\langle#1,#2\right\rangle}
\newcommand{\lp}{\left(}
\newcommand{\rp}{\right)}
\newcommand{\lb}{\left[}
\newcommand{\rb}{\right]}
\newcommand{\lc}{\left\{}
\newcommand{\rc}{\right\}}
\newcommand{\lab}{\left|}
\newcommand{\rab}{\right|}

\newcommand{\Lap}{\Delta}

\newcommand{\Hess}{\nabla^{2}}
\DeclareMathOperator{\dist}{dist}
\DeclareMathOperator{\diam}{diam}

\newcommand{\sO}{\mathcal{O}}
\newcommand{\sH}{\mathcal{H}}
\DeclareMathOperator{\vol}{vol}

\DeclareMathOperator{\Var}{Var}
\newcommand{\E}{\mathbb{E}}
\newcommand{\Emu}[1]{\E^{\mu}\lb #1 \rb}
\newcommand{\Varmu}[1]{\Var^{\mu}\lb #1 \rb}

\newcommand{\GE}{\Gamma_{\epsilon}}

\DeclareMathOperator{\Cut}{Cut}
\DeclareMathOperator{\sing}{sing}

\DeclareMathOperator{\Sym}{Sym}
\DeclareMathOperator{\Area}{Area}

\DeclareMathOperator{\Lip}{Lip}

\renewcommand{\d}[1]{\,d#1}
\renewcommand{\phi}{\varphi}

\newcommand{\St}{S_{1/2}(t_i)}
\newcommand{\SetChar}{\boldsymbol 1}
\newcommand{\tz}{\tilde{z}}

\begin{document}

\title[heat kernel at cut locus]{The small time asymptotics of the heat kernel \\ at the cut locus}
\author{Robert Neel}
\address{Department of Mathematics, Columbia University, New York, NY}
\begin{abstract}We study the small time asymptotics of the gradient and Hessian of the logarithm of the heat kernel at the cut locus, giving, in principle, complete expansions for both quantities.  We relate the leading terms of the expansions to the structure of the cut locus, especially to conjugacy, and we provide a probabilistic interpretation in terms of the Brownian bridge.  In particular, we show that the cut locus is the set of points where the Hessian blows up faster than $1/t$.  We also study the distributional asymptotics and use them to compute the distributional Hessian of the energy function (that is, one-half the distance function squared). \end{abstract}
\thanks{The author gratefully acknowledges support from an NSF Graduate Research
  Fellowship and a Clay Liftoff Fellowship.  This paper is partially based on the author's Ph.D.\ thesis.}
\email{neel@math.columbia.edu}
\date{May 29, 2006}

% AMS notes
% 58J35 Heat and other parabolic equation methods
% 58J65 Diffusion processes and stochastic analysis on manifolds
% 53C22 Geodesics 

\maketitle

\section{Introduction}

Let $M$ be a compact, connected, smooth Riemannian manifold of
dimension $n$.  For any point $x\in M$, we use $\Cut(x)\subset M$ to
denote the cut locus of $x$.  In particular, since $M$ is compact,
$\Cut(x)$ will be nonempty for every $x$.  The
Riemannian metric induces a distance function $\dist(x,y)$.  We will
also need to consider the energy function,
$E(x,y)=\frac{1}{2}\dist(x,y)^2$.  Let $\Lap$ be the Laplace-Beltrami operator on $M$, that is, if $x_1,\ldots,x_n$ are normal coordinates centered at a point $p$, then
$\Lap = \sum_{i=1}^{n} \frac{\partial^2}{\partial x_{i}^{2}}$ at $p$.
The heat kernel $p_t(x,y)$ is the fundamental solution to the heat
equation $\partial_t u(t,x) = \frac{1}{2}\Lap u(t,x)$.

A well-known result of Varadhan states that $t\log p_{t}(x,y)$ converges to $-E(x,y)$ as $t\da 0$ uniformly on all of $M$.  Motivated by this result, we define
\[
E_t(x,y) = -t \log p_t(x,y)
\]
and observe that $E_t(x,y)\ra E(x,y)$ uniformly.  Malliavin and Stroock \cite{SAndM} have shown, using pathspace methods, that away from the cut locus, spatial derivatives of $E_t(x,y)$ commute with taking the limit as $t\da 0$ (for an analytic proof, see~\cite{BGV}).  Clearly, the lack of differentiability of $E(x,y)$ at the cut locus means that something else must be occurring there (see Bishop \cite{Bishop} for a brief discussion of the smoothness properties of the distance function at the cut locus).  Indeed, in the same paper, Malliavin and Stroock use pathspace
integration to show that, if the set of minimal geodesics connecting
$x$ and $y$ is sufficiently ``nice,'' then $\Hess E_t(x,y)$
is asymptotic to $-1/t$ times the variance of a random variable on path space as $t\da 0$.  Unfortunately, their analysis is too complicated to obtain more detailed information.

In the present paper, we develop an analogous, but purely finite-dimensional, approach which allows a much more detailed analysis of the small-time asymptotics of the gradient and Hessian of $E_t(x,y)$.   In particular, we will show how complete asymptotic expansions of the gradient and the Hessian of $E_t(x,y)$ can be represented as an integrals over the set of midpoints of minimal geodesics from $x$ to $y$.  These general expansions are given in Theorems~\ref{THM:Grad} and ~\ref{THM:Hess} respectively.  The leading term in each of these expansions, given in Equation~\eqref{Eqn:Leading}, is fairly accessible both to analysis and to a probabilistic interpretation, and we will show how the small time behavior of $\Hess E_t(x,y)$ reflects the structure of minimal geodesics from $x$ to $y$.  In particular, we will show (see Theorem~\ref{THM:Char} below) that $\Hess E_t(x,y)$, as a quadratic form on $T_yM$, is unbounded as $t\da 0$ if and only if $y\in\Cut(x)$.  Finally, in Theorems~\ref{THM:rho} and~\ref{THM:CPR}, we show how the asymptotic expansion of $\Hess E_t(x,y)$ can be used to determine the distributional Hessian of $E(x,y)$.  

We note that our methods are an extension of a procedure due to Molchanov \cite{Molchanov} of representing the heat kernel itself as an integral over the midpoints of minimal geodesics, allowing him to determine the rate of decay of the heat kernel at points $y\in \Cut(x)$ in a variety of cases.

A survey of the results in the present paper can be found in~\cite{MeStroock}.

I would like to thank my advisor, Dan Stroock, for his invaluable suggestions throughout the course of this work.  I also thank David Jerison and Joe Harris for helpful discussions.

\section{The representation as Laplace integrals}

Our main tool will be a pair of formulas which express the gradient and the Hessian of $E_t(x,y)$ as integrals over the set of midpoints of minimal geodesics from $x$ to $y$, the asymptotics of which are amenable to study.  The present section is devoted to the derivation of these formulas.

\subsection{Preliminary results}  We will need various facts about the small time asymptotics of the heat kernel, which we present here.

Let $C_{M} \subset M \times M$
be the set of pairs of points $(x,y)$ such that $y\in \Cut(x)$.  Away from the cut locus, we have the following asymptotic expansion of the heat
kernel, due to Minakshisundaram and Pleijel \cite{MP} (see \cite{Chavel} for a more modern development).
\begin{THM}
\label{THM:Pleijel}
Let $M$ be a smooth, complete Riemannian manifold of dimension $n$.
Then there are smooth functions $H_{i}(x,y)$ defined on $(M\times M)\backslash C_M$ such that
the asymptotic expansion
\[
p_{t}(x,y) \sim \lp\frac{1}{2\pi t}\rp^{n/2} e^{-E(x,y)/t}
\sum_{i=0}^{\infty} H_{i}(x,y)t^{i}
\]
holds uniformly as $t\da 0$ on compact subsets of $(M\times
M)\backslash C_M$.
Further, if $y=\exp_{x}(Y)$, then $H_{0}(x,y)$ is given by the
reciprocal of the square root of the Jacobian of $\exp_{x}$ at $Y$.
\end{THM}
For future use, let $k(t,x,y)=(2\pi t)^{n/2}e^{E(x,y)/t} p_t(x,y)$ be
defined away from the cut locus, so
that $k(t,x,y)\sim \sum_{i=0}^{\infty} H_{i}(x,y)t^{i}$.  Recall the result of Varadhan mentioned in introduction; namely that $E_t$ converges to $E$ uniformly on all of $M$.  We will rewrite this as
\be
\label{Eqn:Varadhan}
p_{t}(x,y) =\exp\lb \frac{-E(x,y)+\delta(t,x,y)}{t} \rb
\ee
where $\delta(t,x,y)$ is some function which
goes to 0 uniformly in $t$ on all of $M$.  Incidentally, one cannot hope
to replace $\delta(t,x,y)$ with a power series expansion.  Even in
the simplest case of the heat kernel on $\mS^1$, we see that
$\delta(t,x,y)$ fails to be $O(t)$ for any $x$ and $y$.

Having summarized the small-time asymptotics of the heat kernel itself
both away from and on the cut locus, we now turn to the log derivatives.
As mentioned, Malliavin and Stroock \cite{SAndM} (for $m=1, 2$) and Stroock and Turetsky \cite{SAndT2} (for $m>2$) have proved that
\be
\label{Eqn:SAndM}
\nabla^{m} E_{t}(x,y) \ra \nabla^{m} E(x,y)
\ee
uniformly on compact subsets of $M\backslash \Cut(x)$, where
$\nabla^{m}$ is the mth covariant derivative and all derivatives are
taken in the $y$ variable.
Next, we will need better control over the convergence of $\nabla^m E_{t}(x,y)$ away from the cut locus.
\begin{Lemma}
\label{Lemma:Gi}
Let $M$ be a smooth, complete Riemannian manifold of dimension $n$.
Then there are smooth functions $G_{i}(x,y)$ defined on $(M\times M)\backslash C_M$ such that, for any positive integer $m$,
the asymptotic expansion
\[
t\nabla^m \log p_{t}(x,y) \sim \sum_{i=0}^{\infty} \nabla^m G_{i}(x,y)t^{i}
\]
holds uniformly as $t\da 0$ on compact subsets of $(M\times
M)\backslash C_M$.
Further, the $G_i$ are given in terms of the $H_i$ by taking the log derivatives of the 
Minakshisundaram-Pleijel expansion of Theorem \ref{THM:Pleijel}, and in particular, $G_0(x,y)=-E(x,y)$.
\end{Lemma}

This lemma is a direct consequence of the fact that the Minakshisundaram-Pleijel expansion can be differentiated.  For a probabilistic proof of this fact see~\cite{BenArous}; for an analytic proof see~\cite{BGV}.  Note that we get a power series expansion only if we take at least one spatial derivative.   Again for future use, let $l(t,x,y,A)=t\nabla_A\log p_{t}(x,y)$ be defined
away from the cut locus, so that it has the above expansion.

Finally, we have a result of
Stroock and Turetsky~\cite{SAndT} which states that
\be
\label{Eqn:SAndT}
\left|\nabla^{m} p_{t}(x,y)\right|\leq D_{m}\lb\frac{\dist(x,y)}{t}
+\frac{1}{\sqrt{t}}\rb^{m}
p_{t}(x,y)
\ee
where the $D_m$ are some constants depending only on $M$.  The point is
that this estimate is valid even when $y\in\Cut(x)$.  Note that the
results of Stroock and Malliavin mentioned above for $m=2$, and the results of Stroock and Turetsky for $m>2$, show that 
the power of $t$ in Equation~\eqref{Eqn:SAndT} is sharp.

\subsection{The gradient}

In this section, we will prove a theorem describing the asymptotic expansion of the gradient of $E_t(x,y)$ which is valid everywhere on $M$.  We begin by introducing some notation.  Fix any two distinct points $x$ and $y$ on $M$.  Let $\Gamma$ be the set
of midpoints of minimal geodesics from $x$ to $y$ (for example, in the case $M=\mS^n$ with $x$ and $y$ the north and south poles, $\Gamma$ is the equator).  By compactness, there exists
$\epsilon>0$ such that the $\epsilon$-neighborhood of $\Gamma$,
denoted $\Gamma_{\epsilon}$, is strictly positive distance from $x$, $y$, and
both of their cut loci.  Also, we define the \emph{hinged energy function} $h_{x,y}(z)=E(x,z)+E(y,z)$.  Note that $h_{x,y}$ obtains its minimum precisely on $\Gamma$ and that this minimum is equal to $E(x,y)/2$.

\begin{THM}
\label{THM:Grad}
Let $M$ be a smooth, compact, connected Riemannian manifold.  Choose any two distinct points $x$ and $y$.  Then there exist positive constants $C$ and
$\lambda$ such that, for any $A\in T_{y}M$, we have
\[ 
\nabla_{A} E_{2t}(x,y)=
-\frac{\int_{\Gamma_{\epsilon}}
l(t,z,y,A)
\exp\lb-\frac{h_{x,y}(z)}{t}\rb k(t,x,z)k(t,y,z) \d{z} }
{\int_{\Gamma_{\epsilon}} \exp\lb-\frac{h_{x,y}(z)}{t}\rb
k(t,x,z)k(t,y,z) \d{z} } +\hat{e}(t,x,y)
\]
where $l(t,x,y,A)$ and $k(t,x,y)$ are as above and $|\hat{e}(t,x,y)|\leq
C\exp\lp-\lambda/t\rp$.
\end{THM}
\emph{Proof.}  Choose and fix some $(x,y)\in M\times M$.  Let $\Gamma$ and $\Gamma_{\epsilon}$ be as above with $\epsilon$ small enough so that $\Gamma_{\epsilon}$ is a strictly positive
distance from $\Cut(x)$ and $\Cut(y)$.  The Chapman-Kolmogorov equation
gives
\[
p_{2t}(x,y) = \int_{\Gamma_{\epsilon}}p_t(x,z)p_t(z,y)\d{z}
  + \int_{M\backslash\Gamma_{\epsilon}}p_t(x,z)p_t(z,y)\d{z} .
\]
Using Theorem~\ref{THM:Pleijel} and Equation~\eqref{Eqn:Varadhan}, we have
\begin{multline*}
p_{2t}(x,y) = \int_{\Gamma_{\epsilon}}\lp\frac{1}{2\pi t}\rp^{n}
\exp\lb-\frac{h_{x,y}(z)}{t}\rb k(t,x,z)k(t,y,z) \d{z} \\
+ \int_{M\backslash\Gamma_{\epsilon}} \exp\lb-\frac{1}{t}
\lp h_{x,y}(z)+\delta(t,x,y) \rp\rb \d{z} .
\end{multline*}
Observe that $h_{x,y}(z)$ is continuous and achieves its minimum
precisely on $\Gamma$.  Hence its minimum on
$M\backslash\Gamma_{\epsilon}$ is strictly greater than its minimum on
$\Gamma_{\epsilon/2}$.  It follows that there
exist positive $\lambda$ and $C$ such that
\be
\label{Eqn:p2t}
p_{2t}(x,y) = \lb1+e(t,x,y)\rb\int_{\Gamma_{\epsilon}}\lp\frac{1}{2\pi T}\rp^{n}
\exp\lb-\frac{h_{x,y}(z)}{t}\rb k(t,x,z)k(t,y,z) \d{z}
\ee
where $|e(t,x,y)|\leq C\exp\lp-\lambda/t\rp$ (see Lemma 5.3.1 of Hsu's
book~\cite{Hsu}, where
he gives essentially this result with a more detailed proof).

Again use the Chapman-Kolmogorov equation to write
\[
p_{2t}(x,y)=\int_{M}p_{t}(x,z)p_{t}(z,y) \d{z}.
\]
Taking derivatives (all derivatives are on the $y$ variable) gives
\[
\nabla_{A}p_{2t}(x,y)=\int_{M}p_{t}(x,z)\lb\nabla_{A}p_{t}(z,y)\rb
\d{z} .
\]

We divide the manifold into
three regions.  Let $\epsilon>0$, $\Gamma$, and $\Gamma_{\epsilon}$ be
as before.  Let $C_{\epsilon}$ be
an $\epsilon$-neighborhood around the cut locus of $x$.  We
now demand that $\epsilon$ also be small enough that these sets are a
strictly positive distance from one another.  Finally, let
$R_{\epsilon}=M\backslash(C_{\epsilon}\cup\Gamma_{\epsilon})$ be the rest of $M$.
Then
\begin{multline*}
\nabla_{A}p_{2t}(x,y) = \int_{\Gamma_{\epsilon}}p_{t}(x,z)
\lb\nabla_{A}p_{t}(z,y)\rb \d{z}
+ \int_{R_{\epsilon}}p_{t}(x,z)\lb\nabla_{A}p_{t}(z,y)\rb
\d{z} \\
 + \int_{C_{\epsilon}}p_{t}(x,z)\lb\nabla_{A}p_{t}(z,y)\rb
\d{z} .
\end{multline*}
On $\Gamma_{\epsilon}$ and $R_{\epsilon}$ we are a strictly positive
distance from the cut locus.  Also, note that
\[
\nabla_{A}p_{t}(z,y)=p_{t}(z,y) \nabla_{A}\log p_{t}(z,y)
\]
Then we have
\begin{multline*}
\nabla_{A}p_{2t}(x,y)  =
\int_{\Gamma_{\epsilon}\cup R_{\epsilon}}
\frac{1}{t}l(t,z,y) \lp\frac{1}{2\pi t}\rp^{n}
 \exp\lb-\frac{h_{x,y}(z)}{t}\rb k(t,x,z)k(t,y,z)\d{z} \\ 
  +\int_{C_{\epsilon}}p_{t}(x,z)\lb\nabla_{A}p_{t}(z,y)\rb \d{z} .
\end{multline*}

So now the problem is to control the last integral.  For this we
use Equations~\eqref{Eqn:Varadhan} and~\eqref{Eqn:SAndT}.
This gives us the bound
\begin{multline*}
\left|\int_{C_{\epsilon}}p_{t}(x,z)\lb\nabla_{A}p_{t}(z,y)\rb
\d{z}\right|  \leq  \\
 D_1\lb\frac{\diam(M)}{t} +\frac{1}{\sqrt{t}}\rb
\int_{C_{\epsilon}}
\exp\lb-\frac{1}{t}
\lp h_{x,y}(z) +\delta(t,x,y)\rp\rb \d{z} .
\end{multline*}
To get the log gradient, we need to divide through by $p_{2t}(x,y)$,
which we write as an integral using Equation~\eqref{Eqn:p2t}.
We now claim that the terms involving $R_{\epsilon}$ and
$C_{\epsilon}$ vanish exponentially fast.  Indeed, we've already seen this
during the derivation of Equation~\eqref{Eqn:p2t}.  Thus
we have
\begin{multline*}
\nabla_{A} \log p_{2t}(x,y)= \\
\frac{1}{1+e(t,x,y)} \left\{
\frac{1}{t}
\frac{\int_{\Gamma_{\epsilon}}
l(t,z,y)
\exp\lb-\frac{h_{x,y}(z)}{t}\rb k(t,x,z)k(t,y,z) \d{z} }
{\int_{\Gamma_{\epsilon}} \exp\lb-\frac{h_{x,y}(z)}{t}\rb
k(t,x,z)k(t,y,z) \d{z} } +\hat{e}(t,x,y) \right\}
\end{multline*}
where $e(t,x,y)$ and $\hat{e}(t,x,y)$ are both bounded in absolute
value by $C\exp\lp-\lambda/T\rp$, for some positive $C$ and $\lambda$
(perhaps different from above).  Further, we have that $1/(1+e(t,x,y))
=1+O(e(t,x,y))$.  Hence, by adjusting $\hat{e}(t,x,y)$, $C$, and
$\lambda$, we can get rid of $e(T,x,y)$.  Recalling the definition of $E_t(x,y)$, we have proved the theorem.  $\Box$

\subsection{The Hessian}

Here we develop the analogous formula for the Hessian.  By polarization, it is sufficient to consider $\Hess_{A,A}E_{2t}(x,y)$.
\begin{THM}
\label{THM:Hess}
Let $M$ be a smooth, compact, connected Riemannian manifold.  Choose
any two distinct points $x$ and $y$.  Then there exist positive constants $C$ and
$\lambda$ (possibly different from the constants in Theorem~\ref{THM:Grad}) such that, for any $A\in T_{y}M$, we have
\begin{multline*}
\Hess_{A,A}E_{2t}(x,y) = 
-\frac{1}{t} \left\{
\frac{\int_{\Gamma_{\epsilon}}
\lp l(t,z,y,A)\rp^{2}
\exp\lb-\frac{h_{x,y}(z)}{t}\rb
k(t,x,z)k(t,y,z) \d{z} }
{\int_{\Gamma_{\epsilon}}
\exp\lb-\frac{h_{x,y}(z)}{t}\rb
k(t,x,z)k(t,y,z) \d{z}} \right.  \\
  \left.  -\lb
\frac{\int_{\Gamma_{\epsilon}}
l(t,z,y,A)
\exp\lb-\frac{h_{x,y}(z)}{t}\rb
k(t,x,z)k(t,y,z) \d{z} }
{\int_{\Gamma_{\epsilon}}
\exp\lb-\frac{h_{x,y}(z)}{t}\rb
k(t,x,z)k(t,y,z) \d{z}} \rb^{2} \right\} \\
-\frac{\int_{\Gamma_{\epsilon}}
\nabla_{A}l(t,z,y,A)
\exp\lb-\frac{h_{x,y}(z)}{t}\rb
k(t,x,z)k(t,y,z) \d{z} }
{\int_{\Gamma_{\epsilon}}
\exp\lb-\frac{h_{x,y}(z)}{t}\rb
k(t,x,z)k(t,y,z) \d{z} } +\hat{e}(T,x,y)
\end{multline*}
where $l(t,x,y,A)$ and $k(t,x,y)$ are as above and $|\hat{e}(t,x,y)|\leq
C\exp\lp-\lambda/t\rp$.
\end{THM}
\emph{Proof.}
We begin by observing that
\be
\label{Eqn:log}
\Hess_{A,A} \log f 
 =  \frac{\Hess_{A,A}f}{f}- (\nabla_{A}\log f)^{2} .
\ee
Since we have computed the log gradient of the heat kernel in our proof of
Theorem~\ref{THM:Grad}, the only thing remaining is for us to compute
the Hessian of the heat kernel.  Again we start with the
Chapman-Kolmogorov equation and differentiate under the integral sign to get
\[
\begin{split}
\Hess_{A,A}p_{2t}(x,y)
 & =
 \int_{M} p_{t}(x,z) \Hess_{A,A} p_{t}(z,y) \d{z} \\
 & = \int_{\Gamma_{\epsilon}\cup R_{\epsilon}} \lb \Hess_{A,A} \log p_{t}(z,y)
 + \lp \log p_{t}(z,y)\rp^{2}\rb p_{t}(z,y) p_{t}(z,x) \d{z} \\
& \quad + \int_{C_{\epsilon}} p_{t}(z,x) \Hess_{A,A} p_{t}(z,y) \d{z} \\
 &= \int_{\Gamma_{\epsilon}\cup R_{\epsilon}}
\lb\frac{1}{t}\nabla_{A}l(t,z,y,A)+
\frac{1}{t^2}\lp l(t,z,y,A)\rp^{2} \rb \times \\
 & \quad \lp\frac{1}{2\pi t}\rp^{n} \exp\lb-\frac{h_{x,y}(z)}{t} \rb
 k(t,x,z)k(t,y,z) \d{z} \\
 & \quad + \int_{C_{\epsilon}}
 p_{t}(z,x) \Hess_{A,A} p_{t}(z,y) \d{z} .
\end{split}
\]

We wish to divide both sides by
$p_{2t}(x,y)$, since that's what appears in the expansion of the log
Hessian.  We can use Equation~\eqref{Eqn:SAndT}
to control the integral over $C_{\epsilon}$ (in particular, it
decays exponentially), and the integral over
$R_{\epsilon}$ also decays exponentially, for the same reasons as before.
Thus (refer to Theorem~\ref{THM:Grad})
\begin{multline*}
\frac{\Hess_{A,A}p_{2t}(x,y)}{p_{2t}(x,y)} = \\
\frac{1}{1+e(t,x,y)} \left\{
\frac{1}{t}\frac{\int_{\Gamma_{\epsilon}}
\nabla_{A}l(t,z,y,A) \exp\lb-\frac{h_{x,y}(z)}{t}\rb
k(t,x,z)k(t,y,z) \d{z}}
{\int_{\Gamma_{\epsilon}}
\exp\lb-\frac{h_{x,y}(z)}{t}\rb
k(t,x,z)k(t,y,z) \d{z}} \right. \\
\left. +\frac{1}{t^2}\frac{\int_{\Gamma_{\epsilon}}
\lb l(t,z,y,A) \rb^{2} \exp\lb-\frac{h_{x,y}(z)}{t}\rb
k(t,x,z)k(t,y,z) \d{z}}
{\int_{\Gamma_{\epsilon}}
\exp\lb-\frac{h_{x,y}(z)}{t}\rb
k(t,x,z)k(t,y,z) \d{z}} +\hat{e}(t,x,y) \right\}
\end{multline*}
where we change $C$ and $\lambda$ as necessary.

We now plug all of our results from above into
Equation~\eqref{Eqn:log}.  The theorem then follows.  $\Box$

These theorems work well when we wish to compute the asymptotics of the log
gradient or log
Hessian with respect to fixed $x$ and $y$.  However, if we wish to
study how the asymptotics change as $y$ moves, say into the cut locus,
then they won't be of much use.  This is because $\Gamma_{\epsilon}$
and $\lambda$ can change discontinuously in $y$, which we can see just
by looking at $\mS^{2}$ and letting $y$ move into $\Cut(x)$.  Fortunately,
our derivation of the theorem makes it clear how to solve this
problem.
Let $\sO$ be the union of the sets $\Gamma$ associated to every $y\in\Cut(x)$.  Then $\sO$ is still a uniformly positive distance from $\Cut(x)$.  Let $\sO_{\epsilon}$ be an
$\epsilon$-neighborhood, chosen small enough to still be a uniformly positive
distance from $\Cut(x)$.  Then if we choose $y\in B_{\epsilon}(\Cut(x))$, everything we've done above works with
$\Gamma$ and $\Gamma_{\epsilon}$ replaced by $\sO$ and
$\sO_{\epsilon}$.  This will allow us
to reduce all questions of what happens when $y$ is moved to
studying how our integral operators change.  Said informally, all we've done is
take all the parts of $R$ which will become relevant as we move $y$
and make them part of $\sO$, our new region of interest.  Thus
all of the important behavior takes places in $\sO$.  This
modification is worth formalizing.
\begin{cor}
\label{TheCor}
Let $M$ and $x$ be as above.  Then for
any $y\in B_{\epsilon}(\Cut(x))$, the expansions in
Theorem~\ref{THM:Grad} and Theorem~\ref{THM:Hess}
hold with $\Gamma_{\epsilon}$ replaced throughout by
$\sO_{\epsilon}$.
\end{cor}

\section{The leading terms}

Theorems~\ref{THM:Grad} and~\ref{THM:Hess} give, in principle, the
complete asymptotic expansions of the gradient and Hessian of $E_t(x,y)$, up to terms which vanish faster than any power of $t$.  However, if we restrict our attention to the first few terms of these
expansions, the formulas simplify considerably.

\subsection{Formulation in terms of expectation and variance}

Considering the leading terms for the gradient and Hessian of $E_t(x,y)$, we have
\begin{eqnarray*}
\lefteqn{\nabla_{A} E_{t}(x,y)=} & & \\
 & &
2\left\{
\frac{\int_{\Gamma_{\epsilon}}
\nabla_{A} E(z,y) \exp\lb-\frac{2}{t}
h_{x,y}(z)\rb H_{0}(x,z) H_{0}(y,z) \d{z}}
{\int_{\Gamma_{\epsilon}} \exp\lb-\frac{2}{t}
h_{x,y}(z)\rb H_{0}(x,z) H_{0}(y,z) \d{z} } \right\} +O(t) .
\end{eqnarray*}
and
\begin{eqnarray*}
\lefteqn{\Hess_{A,A}E_{t}(x,y) = } \\
& &
-\frac{4}{t} \left\{ 
\frac{\int_{\Gamma_{\epsilon}} \lp\nabla_{A}E(z,y)\rp^{2}
\exp\lb-\frac{2}{t} h_{x,y}(z)\rb H_{0}(x,z)H_{0}(y,z) \d{z}}
{\int_{\Gamma_{\epsilon}} \exp\lb-\frac{2}{t}
h_{x,y}(z)\rb H_{0}(x,z)H_{0}(y,z)\d{z}} \right. \\
 & &
\mbox{} \left. -\lb\frac{\int_{\Gamma_{\epsilon}} \nabla_{A}E(z,y)
\exp\lb-\frac{2}{t} h_{x,y}(z)\rb H_{0}(x,z)H_{0}(y,z) \d{z}}
{\int_{\Gamma_{\epsilon}} \exp\lb-\frac{2}{t}
h_{x,y}(z)\rb H_{0}(x,z)H_{0}(y,z)\d{z}} \rb^{2}
\right\} +O(1) .
\end{eqnarray*}
Note that we're now looking at $E_t(x,y)$ at time $t$ rather than
at time $2t$.

To begin, we can give a probabilistic interpretation of the constant term
of the gradient and the $1/t$ term of the Hessian.  Let
\be
\label{Eqn:mut}
\begin{aligned}
\mu_t(dz)&=\frac{\SetChar_{\Gamma _\epsilon }(z)}
{Z_t}H(x,z)H(y,z)\exp\left(-\frac{2h_{x,y}(z)}t\right)\,dz\\
\text{where}\quad Z_t&=\int_{\Gamma _\epsilon
}H(x,z)H(y,z)\exp\left(-\frac{2h_{x,y}(z)}t\right)\,dz.
\end{aligned}
\ee
Then $\mu_t$ is a probability measure supported on $\GE$, and the above becomes
\[
\begin{aligned}\label{Eqn:Leading}
\nabla_{A} E_{t}(x,y) &=2 \E^{\mu_t}\lb\nabla_{A}E(z,y)\rb +O(t) \\
\text{and}\quad  \Hess_{A,A}E_t(x,y) &= -\frac{4}{t} \Var^{\mu_t}\lb{\nabla_{A}E(z,y)}\rb
+O(1) .
\end{aligned}
\]
Since $\Gamma_{\epsilon}$ is compact, so is the space of
probability measures supported on $\Gamma_{\epsilon}$ (in the weak
topology).  In particular, if we take any sequence of times decreasing
to 0, then it will have a subsequence $t_{i}$ such that
$\mu_{t_{i}}$ converges to some limit probability measure $\mu$
supported on $\Gamma$.  Let $M_0$ be the set of all such limit
measures.  For any $\mu\in M_0$, we will say that $t_i$ is an associated sequence of times if $\mu_{t_i}$ converges (weakly) to $\mu$.

We know that
\[
\begin{split}
\nabla_{A}E(z,y) &= \dist(z,y)\ip{A}{\nabla\dist(z,y)} \\
&= \frac{1}{2}\dist(x,y)\ip{A}{Y(z)}
\end{split}
\]
for $z\in \Gamma$, where $Y(z)$ is
the (unique) unit vector at $y$ such that
\[
\exp_{y}\lb -Y(z)\dist(z,y)\rb=z .
\]
Further, if
we let $\theta_A(z)$ be the angle between $A$ and $Y(z)$, then we can write $\ip{A}{Y(z)}$
as $|A|\cos\theta_{A}(z)$.  Thus, if we choose any $\mu\in M_0$ and let $t_i$ be an associated sequence of times, we have
\[
\begin{split}
\lim_{i\ra \infty} \nabla_{A} E_{t_i}(x,y) &= |A|\dist(x,y)\Emu{\cos\theta_{A}(z)} \\
\text{and}\quad \lim_{i\ra \infty}t_{i} \Hess_{A,A} E_{t_i}(x,y) &= -|A|^2\dist(x,y)^{2} \Varmu{\cos\theta_{A}(z)} .
\end{split}
\]

Some elementary facts about the log
gradient and log Hessian follow immediately.  If we homothetically scale $M$ by a factor of
$a>0$, then $\lim_{i\ra \infty}\nabla_{A} E_{t_{i}}(x,y)$ is multiplied
by $a$ and $\lim_{i\ra \infty}t_{i} \Hess_{A,A} E_{t_{i}}(x,y)$ is
multiplied by $a^2$.  Also, we have the pair of inequalities
\[\begin{split}
&-|A|\dist(x,y)\leq\liminf_{t\da 0} \nabla_{A} E_{t}(x,y) \leq
\limsup_{t\da 0} \nabla_{A} E_{t}(x,y) \leq |A|\dist(x,y) \\
\text{and}\quad
&-|A|^2 \dist(x,y)^2 \leq \liminf_{t\da 0}t \Hess_{A,A} E_{t}(x,y) \leq 
\limsup_{t\da 0}t \Hess_{A,A} E_{t}(x,y) \leq 0 .
\end{split}\]

\subsection{Relation to path-space integration}

There is a one-to-one correspondence between $\Gamma$ and the set of minimal
geodesics from $x$ to $y$.  This suggests that we think of $z\in\Gamma$
as parameterizing these minimal geodesics, of any $\mu\in M_0$ as a measure on
the set of minimal geodesics, and of $\cos\theta_{A}(z)$ as a function
on the minimal geodesics (in particular, $\cos\theta_{A}(z)$ is the cosine of the angle between $A$ and the geodesic corresponding to $z$).  This
viewpoint can be fleshed out by considering the Brownian bridge.  (Intuitively, the Brownian bridge from $x$ to $y$ at time $t$ is the stochastic process obtained from Brownian motion by conditioning on the particle starting at $x$ and being at $y$ at time $t$; for a more detailed discussion see, for example, Hsu's book~\cite{Hsu}.)  In particular, fix $x$ and $y$, and let $P_t$ be the
measure on path-space corresponding to the Brownian bridge from $x$ to
$y$ at time $t$.  It is well-known that $P_t$ exists and that its finite
marginal distributions are given in terms of the heat kernel.  To be precise, let $X_\tau$ be the map from path-space to $M$ which sends each path to its position at time $\tau\in[0,t]$.  Then for any finite sequence of times $0=\tau_0 < \cdots < \tau_m < \tau_{m+1}=t$, the joint distribution of $X_{\tau_1},\dots,X_{\tau_m}$ under $P_t$ is
\[
\frac{1}{p_t(x,y)}\prod_{i=0}^{m}p_{\tau_{i+1}-\tau_i}(x_i,x_{i+1})
\]
where, of course, $x_0=x$ and $x_{m+1}=y$.
Now consider $X_{t/2}$, and let $\nu_t$ be the
distribution of $X_{t/2}$ under $P_t$.  Then the density of $\nu_t$
with respect to the volume measure on $M$ is
\[
\frac{\d{\nu_t}}{\d{\vol}}(z)=\frac{p_{t/2}(x,z)p_{t/2}(z,y)}{p_{t}(x,y)} .
\]
In order to study the limiting behavior, we integrate against a smooth
test function $\phi(z)$.  Also, we use Equation~\eqref{Eqn:p2t} in the denominator
to get
\[
\E^{\nu_t}\lb \phi(z)\rb =
\frac{\int_{M}\phi(z)p_{t/2}(x,z)p_{t/2}(z,y)\d{z}}
{\lb1+e(t,x,y)\rb\int_{\Gamma_{\epsilon}}\lp\frac{1}{\pi t}\rp^{n}
\exp\lb-\frac{2}{t}h_{x,y}(z)\rb k(t/2,x,z)k(t/2,y,z) \d{z} }.
\]
The contribution from the integral over $M\backslash\GE$ in the
numerator vanishes exponentially, and on $\GE$ we can use the Pleijel
expansion.  Proceeding as before, we conclude that
\[
\E^{\nu_t}\lb \phi(z)\rb =
\frac{\int_{\GE}\phi(z) \exp\lb-\frac{2}{t}h_{x,y}(z)\rb H_{0}(x,z)
  H_{0}(y,z) \d{z}}
{\int_{\GE} \exp\lb-\frac{2}{t}h_{x,y}(z)\rb H_{0}(x,z)
  H_{0}(y,z) \d{z}} +O(t).
\]
It follows that $\mu_{t_i}\ra \mu$ if and only if $\nu_{t_i}\ra \mu$.
So not only can we think of $\mu$ as a measure on the set of minimal
geodesics, but we can also say that it is the natural such measure induced by the Brownian bridge.

The preceding allows us to view our method as a finite dimensional analogue of that used by Stroock and Malliavin, in which they work directly with the measure $P_t$ on pathspace.  According to the heuristics of Feynman-type path integrals, Wiener measure should be thought of as the probability measure on pathspace given by weighting each path $w(\tau)$ by a weight proportional to
\[
\exp\lp -\frac{1}{2t} \int_0^1 |w^{\prime}(\tau)|^2 \d{\tau} \rp
\]
(even though this is not possible in a rigorous sense).  If we consider paths from $x$ to $y$ at time $t$, then as $t\da 0$ the above indicates that $P_t$, which is Wiener measure conditioned to require that the path be at $y$ at time $t$, should be concentrating on paths that minimize energy.  These paths are, of course, the minimal geodesics joining $x$ and $y$.  Thus in the limit, we expect the infinite dimensional pathspace picture to collapse down to a finite dimensional picture on minimal geodesics.  Understanding this collapse on pathspace is somewhat difficult, but the present approach avoids this by working directly on the manifold from the beginning.  

We will study the limiting measure (or measures) almost entirely in the context of geometric analysis, only occasionally remarking on the probabilistic interpretation.  Nonetheless, the probabilistically inclined reader is encouraged to think about the limiting measure $\mu$ (when it exists) as giving the probability that a Brownian particle travelling from $x$ to $y$ ``instantaneously'' does so via a particular minimal geodesic.

\subsection{An example}
One application of Theorem~\ref{THM:Hess} is the explicit computation of the asymptotics of $\Hess E_t(x,y)$ when $y\in\Cut(x)$.  (Of course, one expects such an explicit computation to be possible only in special cases.)  Here we show how this can be done on $\mS^n$. For simplicity, assume that $n\geq 2$.  We note that $y\in\Cut(x)$ if and only if $x$ and $y$ are antipodal points.  Hence, without loss of generality, we take $x$ and $y$ to be the north and south poles, which we denote $N$ and $S$.  Then $\Gamma$
is the equatorial sphere $\mS^{n-1}(1)$.  (By $\mS^n(r)$, we mean the standard $n$-dimensional sphere of radius $r$.)  By symmetry, we see
that $\mu_{t}$ converges to the uniform probability measure on the equatorial sphere (with respect to
the induced volume measure).  Next, let $A$ be any unit vector in
$T_{y}M$ (it doesn't matter which one, again by symmetry).  Decomposing the equatorial sphere into level sets of $\theta_{A}(z)$, we see that the level set
for any $\theta$ is $\mS^{n-2}(\sin\theta)$.

We know that the gradient of $E_t(N,S)$ is zero by symmetry.  Thus we proceed to computing the Hessian.  Let
$\omega_{m}$ denote the volume of the unit sphere of dimension $m$. We have
\[ \begin{split}
\Emu{\cos^{2}\theta_{A}(z)}
&= \frac{1}{\omega_{n-1}}
\int_{\theta=0}^{\pi}\lp\frac{\pi}{2}\cos\theta\rp^{2}
(\sin\theta)^{n-2}\omega_{n-2} \d{\theta} \\
\text{and}\quad
\Emu{\cos\theta_{A}(z)}^{2}
&= \lp\frac{\omega_{n-2}}{\omega_{n-1}}\rp^{2}\frac{\pi^{2}}{4}
\lp \int_{\theta=0}^{\pi}\cos\theta(\sin\theta)^{n-2} \d{\theta} \rp^{2}.
\end{split} \]
This second term vanishes because $\cos\theta$ is anti-symmetric about
$\pi/2$ while $\sin\theta$ is symmetric.  Using this in our formula for the Hessian gives
\begin{eqnarray*}
\lim_{t\da 0}t \lb\Hess_{A,A} E_{t}(N,S)\rb &=&
\frac{\omega_{n-2}\pi^2}{\omega_{n-1}}\int_{\theta=0}^{\pi}(\cos\theta)^{2}(\sin\theta)^{n-2} \d{\theta} \\
 &=& \frac{\pi^2}{n} .
\end{eqnarray*}
We conclude that $\Hess_{A,A} E_t(N,S)\sim -\frac{\pi^2}{nt}|A|^2$
as $t\searrow0$ for any $A\in T_{S}M$ (the above computation assumes that $n\geq 2$, but this formula extends to the case $n=1$ as can be checked easily by hand).

\subsection{The relation to conjugate points}
In the case of $\mS^n$, we were able to determine the limiting measure $\mu$ using only symmetry considerations.  In general, this won't be so easy, and the limiting measures (or measures) will depend on the behavior of $h_{x,y}$ near its minima.  In particular, the limiting measures will be governed by whether or not these minima are degenerate (in the sense of Morse theory, that is, whether or not the Hessian is positive definite) and, if so, how degenerate they are.  However, before discussing the relationship between degeneracy and the asymptotics of $\mu_t$, we wish to relate this degeneracy to the geodesic geometry of the manifold.

We begin by introducing some terminology.  Given a smooth, real-valued function $f$
which is defined in a neighborhood of the origin in $\mathbb R^n$ and a non-negative integer $m$, we will
say {\it $f$ is constant to exactly order $m$ at the origin in
the direction $\xi\in\mS^{n-1}$} if 
\[
(\partial_t)^i \lb f(t\xi) -f(0) \rb_{t=0}
\] 
is zero for $0\le i \leq m$ but is non-zero for $i=m+1$.  In particular, $f$ is constant to exactly order $0$ in the direction $\xi$ if its first derivative in the direction $\xi$ is non-zero, $f$ is constant to exactly order 1 if its first derivative is zero but not its second, and so on.  We will say that $f$ is constant to finite order in the direction $\xi$ if there exists some non-negative integer $m$ with $f$ constant to exactly order $m$, and we will say that $f$ is constant to order at least $m$ in the direction $\xi$ if the derivatives above vanish for $0\leq i\leq m$.

Now, let $\gamma$ be a minimal geodesic
from $x$ to $y$, and take
$(r,\theta_1,\ldots,\theta_{n-1})$ to be a polar coordinate system on $T_xM$ such
that $\gamma(r)=\exp_x(r,0,\ldots,0)$ for
$r\in[0,\dist(x,y)]$.  We then say that $\gamma$ is {\it conjugate to exactly
order $m$ in the direction $\xi\in\mS^{n-2}$} if $\theta	\rightsquigarrow
\exp_x(\dist(x,y),\theta)$ is constant to exactly order $m$ in the direction $\xi$ (here $\mS^{n-2}$ is the unit tangent space in the tangent space to $\mS^{n-1}$ at $\theta=0$ and thus the Jacobi field induced by $\xi$ is perpendicular to $\gamma$).
In particular, $\gamma$ is conjugate to exactly order 0 in the direction $\xi$ if it is not conjugate, in the usual sense, in this direction, and $\gamma$ is conjugate to some positive order if it is conjugate in the usual sense, with the order of conjugacy indicating how many derivatives of the exponential map vanish in that direction.  We will use the terms conjugate to finite order and conjugate to order at least $m$ analogously to the above.

The relationship between the
degeneracy of the minima of $h_{x,y}$ and the conjugacy of the corresponding geodesics is contained in the following lemma.

\begin{Lemma}
\label{Lem:Conjugate}
Choose distinct points $x$ and $y$ on $M$ and let $d=\dist(x,y)$.  Let
$(r,\theta_1,\ldots,\theta_{n-1})$ and $\gamma$ be as above, and choose some $\xi\in\mS^{n-2}$.  The coordinates $(r,\theta_1,\ldots,\theta_{n-1})$ on the tangent space induce coordinates in a neighborhood of $\gamma$ under the exponential map, and thus we can also ask to what order $h_{x,y}$ is constant at $\exp(d/2,0,\ldots,0)$ in the direction $\xi$.  Then $h_{x,y}$ is constant to finite order in the direction $\xi$ if and only if $\gamma$ is conjugate to finite order in the direction $\xi$.  In this case, there exists a non-negative integer $m$ such that $h_{x,y}$ is constant to exactly order $2m+1$ at $\exp(d/2,0,\ldots,0)$ in the direction $\xi$ and $\gamma$ is conjugate to exactly order $2m$ in the direction $\xi$.
\end{Lemma}

\emph{Proof.}  
Let $\phi(t)$ be the angle between the geodesic from $\exp(d/2,t\xi)$ to $\exp(d,t\xi)$ and the geodesic from $\exp(d/2,t\xi)$ to $\exp(d,0)=y$.
Let
\[
\rho(t)=\dist(\exp(d,t\xi),y) .
\]
Then, considering the dependence of geodesics on their starting points in the tangent bundle, we see that $\phi(t)$ and $\rho(t)$ are comparable for small $t$.  Thus, $\phi(t)$ is constant to exactly order $l$ at $0$ (in the direction $\partial_t$) if and only $\rho(t)$ is constant to exactly order $l$ at $0$ (in the direction $\partial_t$).  From the definition of $\rho(t)$, it's clear that $\rho(t)$ is constant to exactly order $l$ at $0$ if and only $\gamma$ is conjugate to exactly order $l$ in the direction $\xi$.

The distance from $x$ to $\exp(d/2,t\xi)$ is independent of $t$, and the vector field given by pushing $\xi$ forward by the exponential map is always perpendicular to the geodesic from $x$ to $\exp(d/2,t\xi)$.  Thus, $h_{x,y}(\exp(d/2,t\xi))$ depends only on the distance between $\exp(d/2,t\xi)$ and $y$.  Further, we see that the inner product between the push-forward of $\xi$ and the unit tangent to the geodesic from $\exp(d/2,t\xi)$ to $y$ at $\exp(d/2,t\xi)$ is constant to exactly order $l$ (in the direction $\partial_t$ at $t=0$) if and only if $\phi(t)$ is constant to exactly order $l$. Since this inner product is proportional to the derivative of $h_{x,y}(\exp(d/2,t\xi))$ with respect to $t$, it follows that $\phi(t)$ is constant to exactly order $l$ if and only if $h_{x,y}$ is constant to exactly order $l+1$ in the direction $\xi$.

Combining these facts, we see that $\gamma$ is conjugate to exactly order $l$ in the direction $\xi$ if and only if $h_{x,y}$ is constant to exactly order $l+1$ in the direction $\xi$.  The lemma will be proved once we determine that $h_{x,y}$ can only vanish to odd order.  This, however, is just a restatement of the fact that the first non-zero derivative of $h_{x,y}$ in any direction must be even because $h_{x,y}$ has a local minimum at $\exp(d/2,t\xi)$.  $\Box$

In particular this lemma implies that, if $z\in\Gamma$, then $z$ is a non-degenerate
minimum of $h_{x,y}$ if and only if $x$ and $y$ are conjugate to exactly order $0$ in all directions, which, as mentioned, is the same as saying that $x$ and $y$ are not conjugate along $\gamma$ in the usual sense.  On the other hand, if $z$ is
a degenerate minimum, then $x$ and $y$ are
conjugate (in the usual sense) along the minimal geodesic through $z$, and furthermore, the index and orders of conjugacy\footnote{Some authors use the order of a
conjugate point to denote the dimension of the null space of the
differential of the exponential map.  We will, taking our cue from
Morse theory, instead call that the index of the conjugate point and
reserve the term order to denote the order to which the differential
of the exponential map vanishes in some direction.} can be
determined from the partial derivatives of
$h_{x,y}$.

\subsection{Laplace asymptotics and the resolution of singularities}

Having given a geometric interpretation of the minima of $h_{x,y}$, we now turn to the investigation of the relationship between the minima of $h_{x,y}$ and the limiting measures
$M_0$.  To do so, we will need
to delve into the theory of Laplace asymptotics.  We begin by considering the case when $\Gamma$ consists of finitely many points, say $z_1,\ldots,z_m$.  For ease of notation, we write $g(z) = 2h_{x,y}(z)-E(x,y)$; hence $g$ is non-negative and has zeroes precisely at the $z_i$.  By taking $\epsilon$ small enough, we can ensure that $\GE$ is the union of the disjoint balls $B_\epsilon(z_i)$, where $i$ ranges from $1$ to $m$.  In this case, we see have that
\[
\begin{split}
\E^{\mu_t}\lb\phi\rb &=\frac{1}{Z(t)} \sum_{i=1}^{m} \int_{B_{\epsilon}(z_i)} \phi(z) H_0(x,z)H_0(y,z) e^{-g(z)/t}\d{z} \\
\text{where}\quad Z(t) &=\sum_{i=1}^{m} \int_{B_{\epsilon}(z_i)} H_0(x,z)H_0(y,z) e^{-g(z)/t}\d{z} .
\end{split}
\]
We are thus lead to study integrals of the form
\be
\label{Eqn:LInt}
\int_{B_{\epsilon}(z_i)} \phi(z)e^{-g(z)/t}\d{z}
\ee
as $t\da 0$.

First, suppose that $g$ can be diagonalized at $z_i$, that is, suppose that we can find coordinates $u_1,\ldots,u_n$ around $z_i$ such that
\begin{equation}
\label{Eqn:Diag}
g(u_1,\ldots,u_n)= \sum_{j=1}^{n}u_{j}^{2k_j}
\end{equation}
for some positive integers $k_1\leq\cdots\leq k_n$.  Of course, at a
non-degenerate minimum, the Morse Lemma guarantees the existence of such
coordinates with $k_j=1$ for each $j$, but in general this need not be true.  Under the assumption that \eqref{Eqn:Diag} holds, Estrada and Kanwal \cite{EAndK} give the full asymptotic expansion of \eqref{Eqn:LInt}; namely,
\begin{multline}
\label{Eqn:KEExpansionFull}
\int_{B_{\epsilon}(z_i)} \phi(u_1,\ldots,u_n) \exp\lb -\frac{\sum_{j=1}^{n}u_{j}^{2k_j}}{t}\rb \d{u_1}\cdots \d{u_n} \\ 
\sim \prod_{j=1}^{n} \lb\sum_{i_j=0}^{\infty}
\frac{\Gamma\lp\frac{2i_j+1}{2k_j}\rp}
{(2i_j)!k_j} t^{\frac{2i_j+1}{2k_j}} \lp\frac{\partial}{\partial u_j}\rp^{2i_j}\rb \phi(0,\ldots, 0)
\end{multline}
where derivatives on the right-hand side are to be understood as giving partial differential operators which are then applied to $\phi$ and evaluated at the origin (that is, at $z_i$).  We will be interested mainly in the first term (that is, the one coming from $i_j=0$ for all $j$), in which case Equation~\eqref{Eqn:KEExpansionFull} gives
\begin{equation}
\label{Eqn:KEExpansion}
\int_{B_{\epsilon}(z_i)} \phi(z) e^{-g(z)/t} \d{z} =
  t^{1/2k_{1,i}+\cdots +1/2k_{n,i}} \lb c_i \vol_u(z_i) \phi(z_i) +O\lp t^{1/k_{n,i}}\rp \rb
\end{equation}
where $\vol_u$ is the volume element in the $u$ coordinate chart and
$c_i$ is a positive constant which depends only on $n$ and the $k_{j,i}$'s.

Equation~\eqref{Eqn:KEExpansion} implies that the limit is dominated by those geodesics which are ``the most conjugate,'' in the sense of being conjugate in many directions and/or to high order.  More precisely, assume that $g$ can be diagonalized around each of its minima.  Then to each $z_i$ we can
associate the order of the leading term of the integral over $B_{\epsilon}(z_i)$, which is $l_i=1/2k_{1,i}+\cdots +1/2k_{n,i}$.  Then, as $t\searrow0$, $\mu_t$
converges to a limit $\mu$ which is supported on those $z_i$ with the smallest leading order (that is, $z_i$ with $l_i=\min\{l_1,\ldots,l_m\}$).  Further, the mass
at these points is given by the coefficient of the leading term of the
expansion coming from~\eqref{Eqn:KEExpansion}, normalized to have
total mass one.  In terms of the Brownian bridge, this says that that the particle prefers to travel along the most conjugate geodesics.  More precisely, if we require the particle to travel from $x$ to $y$ ``instantaneously,'' then it travels along the geodesic through $z_i$ with probability $\mu\{z_i\}$.  To be a bit more intuitive, it is not be surprising that a particle under Brownian motion prefers conjugate geodesics, since conjugacy should make a geodesic more forgiving of the white noise which the particle experiences as it tries to follow the geodesic.

In order to extend our analysis beyond the diagonalizable case, we will need to introduce resolutions of singularities.  Let $X$ be a real-analytic manifold of dimension $n$ (we do not require that $X$ be compact), and let $f$ be a real-analytic function on $X$ which is not identically zero.  Then, for our purposes, a resolution of singularities will mean an $n$-dimensional analytic manifold $Y$ and a proper, surjective analytic map $\xi:Y\ra X$ such that $\xi$ possesses the following two properties:
\begin{itemize}
\item At each point of the pre-image of the zero level-set of $f$, there exist local coordinates with respect to which $f\circ \xi$ is a monomial and the Jacobian of $\xi$ is a monomial times a non-zero smooth function;
\item The function $\xi$ is a diffeomorphism outside of the zero level-set of $f$.
\end{itemize}
A famous theorem of Hironaka states that, for $X$ and $f$ as above, a resolution of singularities always exists.

{\bf Remark.}  The original results of Hironaka are both much more general than what we have stated above and formulated in the language of schemes.  A similar definition to the above is given in Arnold, Gusein-Zade, and Varchenko's book~\cite{Arnold}, which is also concerned with the evaluation of Laplace asymptotics.  An ``elementary,'' constructive proof in the case of characteristic zero (which includes the case of real-analytic manifolds) is given by Bierstone and Milman~\cite{BierstoneMilman}.  Finally, a purely analytic statement and proof of a weaker form of the resolution of singularities for real-analytic manifolds (in particular, without the property that $\phi$ is a diffeomorphism off of the zero level-set) is given by Sussmann~\cite{Sussmann}.

As a first application of the resolution of singularities, we can prove the following theorem.

\begin{THM}
\label{THM:Unique}
Let $M$ be a real-analytic, compact Riemannian manifold.  Then for any distinct points $x$ and $y$ of $M$, the corresponding sequence of measures $\mu_t$ (see equation \eqref{Eqn:mut})described above  converges to a unique limit $\mu$ as $t\da 0$.
\end{THM}

\emph{Proof.}  That $M$ is real-analytic implies that the function $g$ defined above is real-analytic in a neighborhood of $\Gamma$, and thus for an appropriate choice of $\epsilon$, on $\GE$.  Let $\phi$ be any smooth function with compact support in $\GE$.  It is enough to show that $\E^{\mu_t}\lb\phi \rb$ converges as $t\da 0$, since the limit of $\mu_t$ is determined by its action on smooth test functions.

Following~\cite{Arnold}, we define an elementary Laplace integral over a bounded, open neighborhood of $0\in \Rn$ as an integral of the form
\[
\int_{U} \exp\lp -x_1^{k_1} \cdots x_n^{k_n}/t\rp \left| x_1^{l_1}\cdots x_n^{l_n} \right| \psi (x_1,\ldots,x_n) \d{x_1}\ldots\d{x_n}
\]
where the $k_i$ are non-negative, even integers, at least one of which is non-zero, the $m_i$ are non-negative integers, and $\psi$ is a smooth function of compact support in $U$ (we call $\psi$ the amplitude function).  By Theorem~7.4 of~\cite{Arnold}, for any elementary Laplace integral there exists a distribution $a$, a positive rational number $\alpha$, and a non-negative integer $m$ such that this integral is asymptotic to $a(\psi)t^{\alpha}|\log t|^{m}$ as $t\da 0$.  Further, if $\psi$ is non-negative on $U$ and positive at the origin, then $a(\psi)$ is positive.

We now apply the resolution of singularities to $g$ and $\GE$ (where we think of $\GE$ as a non-compact real-analytic manifold).  Thus we have a real-analytic manifold $Y$ and a map $\xi:Y\ra \GE$ with the properties listed above.  Because $\Gamma$ (which is the zero level-set of $g$) has measure zero and $\xi$ is a diffeomorphism elsewhere, it follows that an integral over $\GE$ can be written as an integral over $Y$ by pulling back everything on $\GE\setminus\Gamma$ (and ignoring the preimage of $\Gamma$ in $Y$, which has measure zero under the pull-back of the volume form).  In particular, we have that
\[
Z(t)= \int_{Y} (H_0(x,\cdot)\circ\xi) (u) (H_0(y,\cdot)\circ\xi) (u) \exp\lp-g\circ\xi(u)/t\rp \xi^{*}(d\vol)(u) .
\]
Next, observe that restricting the integral to the preimage of the closed $\epsilon/2$ neighborhood of $\Gamma$ only introduces an exponentially small error; and further, the same is true if we integrate over any set intermediate between an $\epsilon/2$ neighborhood and $\GE$.  Since this closed neighborhood is compact, so is its preimage under $\xi$.  Hence it can be covered by finitely many coordinate charts (each of which is also contained in the preimage of $\GE$) of the type described above (namely, $g$ is a monomial and the Jacobian of $\xi$ is a monomial times a smooth, non-zero function).  Using a partition of unity subordinate to this cover, we can write $Z(t)$ as a sum of finitely many, say $l$, elementary Laplace integrals (up to exponentially small error).  In addition, we see that the amplitude function $\psi_i$ in each of these integrals is non-negative and positive at the origin (in the local coordinates) because $H_0$, the volume form, and the partition functions have these properties.  Thus, we see that the integrals over these charts are asymptotic to $a_i(\psi_i)t^{\alpha_i}|\log t|^{m_i}$ respectively, for $i=0,\ldots,l$, where the $a_i(\psi_i)$ are all positive.  Picking out the dominant such term or terms, we conclude that there is some positive real number $c$, some positive rational number $\alpha$, and some non-negative integer $m$ such that
\[
Z(t)\sim c t^{\alpha} |\log t|^m \quad \text{as } t\da 0.
\]

In order to determine the numerator in $\E^{\mu_t}\lb \phi\rb$, we follow the same procedure.  The only difference is that now the amplitude functions in the elementary Laplace integrals also include the pull-back of the test function $\phi$ as a factor.  In particular, we have the same coordinate charts as before, and now our amplitude functions are given by $\tilde{\psi_i}=\psi_i(\phi\circ\xi)$.  We can no longer guarantee that the $\tilde{\psi_i}$ are non-negative and positive at the origins of our charts, and thus we cannot guarantee that $a_i(\tilde{\psi_i})$ is non-zero.  However, note that we need only concern ourselves with the charts where $\alpha_i=\alpha$ and $m_i=m$.  If the sum of the $a_i(\tilde{\phi_i})$ corresponding to these charts is non-zero, then call it $\tilde{c}$ and observe that
\[
\int_{\GE} \phi(z) H_0(x,z)H_0(y,z) e^{-g(z)/t}\d{z} \sim \tilde{c} t^{\alpha} |\log t|^m \quad \text{as } t\da 0.
\]
On the other hand, if the sum of these $a_i(\tilde{\phi_i})$ is zero, we cannot necessarily determine the leading term of the expansion for this integral; however, we can assert that this integral is $o(t^{\alpha} |\log t|^m)$.  In either case, dividing by $Z(t)$ and letting $t\da 0$ shows that $\E^{\mu_t}\lb \phi\rb$ converges to some finite limit.  As mentioned above, this completes the proof. $\Box$

Note that this result says nothing about how to determine the limit $\mu$, or more generally, about how it relates to the geodesic geometry of $M$.  We suspect, although we have not proven, that if $M$ is only assumed to be smooth, the limiting measure need not be unique.

\subsection{Newton polyhedra and evaluation of the limit}

Our use of the resolution of singularities in the last section was non-constructive.  In the present section, we discuss how additional assumptions about $g$ allow one to say much more about the limit measure $\mu$.

We return to the case where $M$ is smooth and $g$ has finitely many minima, which we denote by $z_i$, although now we drop the assumption that $g$ can be diagonalized near the $z_i$.  In a series of papers (for
example~\cite{Vasiliev}) which culminate in the monograph~\cite{Arnold},
Arnold and his school have provided a
fairly complete analysis of the asymptotic expansion of
equation~\eqref{Eqn:LInt} in this case.  We summarize the needed results.  First, we need to assume that $g$ vanishes to finite order at
$z_i$. This is always true in the real-analytic category, but in
the smooth category it need not be the case.   Given that $g$ vanishes to
finite order at $z_i$, we can in fact assume that $g$ is real-analytic in a neighborhood of $g$ by taking an appropriate change of coordinates.  In particular, there exist coordinates $u_j$ around $z_i$ such that $g$ is equal to its Taylor expansion in these coordinates.

Before stating the results, we need to introduce some notation.  Let $\Zed_{+}^n$ be the set of all $n$-tuples of non-negative integers.  Also, we will use the multi-index notation, so that, for $\alpha=(\alpha_1,\ldots,\alpha_n)\in\Zed_{+}^n$ we interpret $u^\alpha$ as the monomial $u_1^{\alpha_1}\cdots u_n^{\alpha_n}$.  Then we have that
\[
g=\sum_{\alpha\in \Zed_{+}^n} g_\alpha u^\alpha
\]
where the $g_\alpha\in\Re$ are the coefficients in the Taylor expansion.  Now let $\Re_{+}^n$ be the set of all $n$-tuples of non-negative reals, and consider $\Zed_{+}^n$ to be a subset of $\Re_{+}^n$ in the natural way.  Then the Newton polytope of $g$ is defined to be the subset of $\Re_{+}^n$ determined by taking the convex hull of $\cup(\alpha+\Re_{+}^n)$ where the union is over all $\alpha\in\Zed_{+}^n$ such that $g_\alpha \neq 0$.  Further, the Newton diagram of $g$, which we will denote by $\Delta(g)$, is defined to be the union of the compact faces of the Newton polytope.  (Strictly speaking, both the Newton polytope and Newton diagram of $g$ also depend on our choice of coordinates $u_i$, but for now we will simply consider the coordinates to be given.)

Let $\gamma$ be a face of $\Delta(g)$.  Then we define $g_\gamma$ to be the power series consisting of the monomials whose multi-indices lie on $\gamma$ with the same coefficients as appear in the power series of $g$.  That is,
\[
g_\gamma = \sum_{\alpha\in \gamma\cap\Zed_{+}^n} g_\alpha u^\alpha .
\]
Note that if $\gamma$ is compact, then $g_\gamma$ is a polynomial.  Similarly, we define the principal part of $g$, denoted $g_\Delta$, by
\[
g_\Delta = \sum_{\alpha\in \Delta(g)\cap\Zed_{+}^n} g_\alpha u^\alpha .
\]
We say that $g$ is non-degenerate if for every compact face $\gamma$ of $\Delta(g)$, the polynomials
\[
\frac{\partial g_\gamma}{\partial x_1},\ldots,\frac{\partial g_\gamma}{\partial x_n}
\]
have no common zeroes in $(\Re\setminus 0)^n$.  We now take a closer look at $\Delta(g)$.  Consider the ray $(r,\ldots,r)$ for $r>0$.  This ray intersects $\Delta(g)$ in exactly one point, say $(r_0,\ldots,r_0)$.  We define $p(\Delta(g))=1/r_0$ to be the remoteness of the Newton diagram of $g$.  Now let $k(\Delta(g))$ be the number of degrees of freedom of the supporting hyperplane of $\Delta(g)$ at the point $(r_0,\ldots,r_0)$; we call $k(\Delta(g))$ the multiplicity of the Newton diagram of $g$.

We are now in a position to state the main result.  Combining Theorems~7.6 and~8.6 of~\cite{Arnold} and Theorem~2.6 of~\cite{Vasiliev}, we have the following.
\begin{THM}
\label{THM:Laplace}
Let $g$ be an analytic function defined in a closed neighborhood $\ol{B_{\epsilon}(0)}\subset \Rn$ of the origin with a unique minimum of finite order at the origin.  Assume that $g(0)=0$.  Let $\phi$ be a smooth function on $\ol{B_{\epsilon}(0)}$ such that $\phi(0)\neq 0$.  Then there exist an integer $m$ between $0$ and $n-1$ inclusive, a positive rational number $\alpha$, and a positive real number $c$, all depending only on $g$, such that
\[
\int_{B_{\epsilon}(0)} \phi(z)e^{-g(z)/t}\d{z} = c\phi(0) t^\alpha |\log t|^m +O(t^\alpha |\log t|^{m-1})
\]
as $t\da 0$.  Consider the Newton diagram of $g$ and its remoteness and multiplicity, as defined above.  Then
\[
t^\alpha |\log t|^m \geq t^{p(\Delta(g))} |\log t|^{k(\Delta(g))}
\]
for small enough $t$.  Finally, if $g$ is non-degenerate, then in fact $\alpha=p(\Delta(g))$ and $m=k(\Delta(g))$.
\end{THM}

While we won't prove this theorem, we will comment on its proof and application.  As might be expected, the resolution of singularities plays an important role.  In fact, the central observation is that in the non-degenerate case one can construct a resolution of singularities from the Newton polytope such that the leading term has the properties given in the theorem.  As such, the above represents a special case of the results of the preceding section.  A more direct approach is taken in \cite{Vasiliev}, in which the lead term of the integral is computed directly under the assumption of non-degeneracy.  The approach taken in \cite{Vasiliev} also proves the one-sided estimate on the lead term given in the theorem.  On the other hand, this approach fails to show that the Laplace integral has an expansion of the desired kind in the degenerate case.

It's clear that Theorem~\ref{THM:Laplace} provides us with more precise information about Laplace integrals when the function $g$ is non-degenerate.  For this reason, it is helpful to know that degeneracy is rare, in some sense.  In particular, Lemma 6.1 of \cite{Arnold} asserts that, given a Newton diagram, the set of degenerate principal parts is a proper semi-algebraic subset of the space of all principal parts corresponding to the given diagram, the complement of which is everywhere dense.  Approaching the issue of degeneracy from a different direction, we could ask to what extent a change of coordinates might help.  Recall that the Newton diagram of $g$ and its associated features, including degeneracy, depend on the chosen coordinates.  It is easy to come up with functions which are degenerate in one coordinate system, but which are non-degenerate in another.  Unfortunately, trying to eliminate degeneracy by changing coordinates does not work in general; there exist functions which are degenerate in any coordinate system.\footnote{The two-dimensional case is almost an exception.  In~\cite{Varchenko}, Varchenko proves that for an analytic function with an isolated minimum, there always exists what he calls \emph{adapted coordinates}, in which the leading power of $t$ is equal to the remoteness of the Newton diagram.}

We now return to the problem at hand, namely determining the limiting behavior of $\mu_t$.  If we assume that $g$ vanishes to finite order at each of the points $z_i$, then we can apply Theorem~\ref{THM:Laplace} to conclude that around each $z_i$ the Laplace integral~\eqref{Eqn:LInt} is asymptotic to a constant times $t^{\alpha_i}|\log t|^{m_i}$.  Since $Z(t)$ is the sum of these terms, we see that $\mu_t$ converges to a limiting measure $\mu$ which is supported on those $z_i$ with dominant leading term.

So far we've been assuming that $\Gamma$ consists of a finite number of points.  In general, though, $\Gamma$ can be quite complex.   While we are far from a result which covers all possibilities, we can generalize the preceding a bit more.  In particular, suppose that $\Gamma$ consists of a finite collection of isolated, smooth submanifolds (possibly with
boundary) $N_1,\ldots,N_m$ of $M$.  
Then the integral over $\Gamma _\epsilon $ can be decomposed into integrals
over the $\epsilon $-neighborhoods $(N_i)_\epsilon $ of the individual
$N_i$.  Further, the integral over each $(N_i)_\epsilon $
can be written as an integral in the normal direction followed by an integral in the tangent direction. Hence, if $g$ vanishes to finite order in the normal directions, we can treat the integral in the normal direction at each $z\in
N_i$ by the above methods.  In
particular, if we assume that for each $N_i$ the asymptotics of the integral in the normal
direction has the same leading term $t^{\alpha_i}|\log t|^{m_i}$ at all points, then the integral over $(N_i)_\epsilon $ will have leading term $t^{\alpha_i}|\log t|^{m_i}$, and
$\mu_t$ will converge to a measure $\mu$ which is supported on the union of
those $N_i$'s with dominant leading term.  Further, on each such $N_i$,
$\mu$ will be absolutely continuous with respect to the induced volume
measure on $N_i$.

While this analysis is fairly general, it certainly does not cover every case.  One can easily come up with cases where, for example, $\Gamma$ has accumulation points or $g$ vanishes to infinite order, and where determining the limiting behavior of $\mu_t$ would be quite difficult, if not impossible.

In a different vein, we can also use Theorem \ref{THM:Laplace} to prove the following lemma, which we will need later.

\begin{Lemma}
\label{Lem:Degenerate}
Let $M$ be a smooth, compact, Riemannian manifold, and let $x$ and $y$ be points in $M$ such that $y\in\Cut(x)$.  Suppose there is a point $w\in\Gamma$ and an associated sequence of times $t_i$ decreasing to zero such that $\mu_{t_i}$ converges to the point mass at $w$ as $i\ra\infty$.  Then $\Hess h_{x,y}(w)$ is degenerate.
\end{Lemma}

\emph{Proof.}  Let $\gamma$ be the minimal geodesic from $x$ to $y$ passing through $w$ (so that $w$ is the midpoint of $\gamma$).
It follows from Lemma \ref{Lem:Conjugate} that $\Hess h_{x,y}(w)$ is degenerate if and only if $\gamma$ is conjugate.  For the remainder of the proof, we assume that $\Hess h_{x,y}(w)$ is non-degenerate.  Hence $\gamma$ is not conjugate, and because $y\in\Cut(x)$, there must be at least one other minimal geodesic from $x$ to $y$.  Let $\hat{\gamma}$ be such a geodesic and let $\hat{w}$ be its midpoint.  Then, because $\Hess h_{x,y}(w)$ is non-degenerate, there are some positive constants $\delta$ and $C$ such that
\[
\begin{split}
\frac{\E^{\mu_t}\lb B_{\delta}(\hat{w}) \rb}{\E^{\mu_t}\lb B_{\delta}(w) \rb} &=
\frac{\int_{B_{\delta}(\hat{w})} H_0(x,z)H_0(z,y) e^{-g(z)/t}\d{z}}{\int_{B_{\delta}(w)} H_0(x,z)H_0(z,y) e^{-g(z)/t}\d{z}} \\
& \geq C \quad\text{for all sufficiently small $t$.}
\end{split}
\]
This follows from the fact that, by multiplying $g$ on $B_{\delta}(\hat{w})$ by some positive constant, $g$ on $B_{\delta}(\hat{w})$ can be made smaller than $g$ on $B_{\delta}(w)$ (where this comparison can be made by introducing normal coordinates on each of these balls), and, by further reducing this constant, any problem arising from the $H_0$ or the volume form can be overcome.  Since $\delta$ can be chosen such that $B_{\delta}(\hat{w})$ and $B_{\delta}(w)$ are disjoint, this inequality shows that no limiting measure can be supported only at $w$.  We have shown that if $\Hess h_{x,y}(w)$ is non-degenerate, then no limit measure can be the point mass at $w$.  This proves the lemma.  $\Box$

\subsection{A characterization of the cut locus}

The results of the previous section show that $\mu$ can be a point mass, even though there are multiple minimal geodesics from $x$ to $y$.  For example, consider the case when there are two minimal geodesics from $x$ to $y$, one of which is conjugate and one of which is not.  Then the limiting measure will be a point mass at the middle of the conjugate geodesic, and the $1/t$ term will vanish.  (In terms of the Brownian bridge, if we require the particle to travel from $x$ to $y$ ``instantaneously,'' then with probability 1 it travels along the conjugate geodesic.)  Thus, it is clear that the $1/t$ term in the asymptotics of $\Hess_{A,A}E_t(x,y)$ is not sufficient to determine when $y\in \Cut(x)$.  On the other hand, we have the following result characterizing $\Cut(x)$ in terms of  the asymptotics of $\Hess E_t(x,y)$.  

\begin{THM}
\label{THM:Char}
Let $M$ be a compact, smooth Riemannian manifold, and let $x$ and $y$
be any two distinct points of $M$.  We have that $y \not\in \Cut(x)$
if and only if
\[
\lim_{t\da 0} \Hess E_{t}(x,y) = \Hess E(x,y)
\]
and $y\in \Cut(x)$ if and only if
\[
\limsup_{t\da 0}  \|\Hess E_{t}(x,y)\|=\infty
\]
where $\|\Hess E_{t}(x,y)\|$ is the operator norm, that is, the max of
$|\Hess_{A,A}E_{t}(x,y)|$ over all $A\in T_{y}M$ with unit length.  Further, if $M$ is real-analytic, we have the stronger result that $y\in \Cut(x)$ if and only if
\[
\lim_{t\da 0}  \|\Hess E_{t}(x,y)\|=\infty .
\]
\end{THM}
\emph{Proof.}
The case where $y \not\in \Cut(x)$ is just a restatement of the result of 
Stroock and Malliavin mentioned in Equation~\eqref{Eqn:SAndM}.  So we consider the case when $y\in \Cut(x)$.  First
suppose that the set of accumulation points of $\mu_t$ contains at
least one measure $\mu$ which is not a point mass.  Let $t_i$ be an
associated sequence of times such that $\mu_{t_i}\ra \mu$.  Because
$\mu$ is not a point mass and the support of $\mu$ is disjoint from
$\Cut(y)$, there is some unit vector $A\in T_yM$ such that $\nabla_{A}E(z,y)$ is not
almost surely $\mu$-constant.  Hence, for this $A$, we
will have a nonzero variance and so $\Hess_{A,A}\log p_{t_i}(x,y)$
will blow up like $1/t^{2}_{i}$.

Next, suppose that all of the accumulation points of $\mu_t$ are point
masses.  This will be the situation for the remainder of the proof.
Take any $w$ such that the point mass at $w$ is an accumulation point of $\mu_t$,
and let $t_i$ be an associated sequence of times (that is, the $t_i$
are a sequence of times decreasing monotonely to zero such that
$\mu_{t_i}$ converges to a point mass at $w$).  Take a small ball
$B_{\delta}(w)$ around $w$, and let
$x_1,\ldots,x_n$ be any smooth coordinates around $w$ defined on this
ball such that $\partial_{x_1},\ldots,\partial_{x_n}$ form an
orthonormal basis for $T_wM$.  Then we wish to consider the function
\[
V(t_i)=\max_{j\in\{1,\ldots,n\}}\Var^{\mu_{t_i}}\lp x_j\rp
\]
where here we consider the $x_j$ to be equal to zero outside of
$B_\delta(w)$.

Define $\rho_{w,\delta}$ to be equal to $2\lb h_{x,y}-h_{x,y}(w)\rb$ on $B_{\delta}(w)$ and zero elsewhere,
and let $\pi_{w,\delta}(s)$ be the measure (with respect to Riemannian volume) of
the subset of $B_{\delta}(w)$ where $\rho_{w,\delta}$ is less than $s$.
Note that $H_0(x,\cdot)H_0(\cdot,y)$ is bounded above and below by positive constants on $\GE$.  In particular, this means that, for the purposes of integration, we can absorb this factor into the volume form.  Doing so greatly simplifies the formulas, and thus we will assume this is the case for the remainder of the proof.  As a result, for the remainder of the proof we let $Z(t)$ be the integral of $\exp\lb -2(h_{x,y}-h_{x,y}(w))/t\rb$ over $\GE$.  (Note that $Z(t)$ is defined in terms of the integral over
all of $\GE$, instead of restricting to $B_{\delta}(w)$.)
Finally, we set
\[
f_{t_i}(x_1,\ldots,x_n) =  \exp\lb -\rho_{w,\delta}(x_1,\ldots,x_n)/t_i\rb/Z(t_i)
\]
on $B_{\delta}(w)$ and equal to zero elsewhere, so that $f_{t_i}$ is the density of $\mu_{t_i}$ with respect to Riemannian volume on
$B_{\delta}(w)$.  (Technically, $f_{t_i}$ also depends on $w$ and $\delta$, but to include this in the notation would be unmanageable.)

We wish to estimate $V(t_i)$ in terms of the measure of the set where
$f_{t_i}$ is at least half of its maximum value,
$f_{t_i}(0,\ldots,0)=1/Z(t_i)$.  Call this set $\St$.
First, we observe that the measure of $\St$ is equal to
$\pi_{w,\delta}(t_i\log2)$.  This follows from noting that
\[
\begin{split}
& \lc(x_1,\ldots,x_n): f_t(x_1,\ldots,x_n)
 \geq \frac{1}{2}f_t(0,\ldots,0)\rc \\
 & \qquad\qquad = \lc(x_1,\ldots,x_n): \frac{\exp\lb
 \rho_{w,\delta}(x_1,\ldots,x_n)\rb}{Z(t)}\geq \frac{1}{2Z(t)} \rc \\
 & \qquad\qquad = \lc (x_1,\ldots,x_n): \rho_{w,\delta}(x_1,\ldots,x_n)\leq t\log2\rc .
\end{split}
\]
Next, we have that
\[
V(t_i) \geq \max_{j\in\{1,\ldots,n\}}\lb \min_{\alpha\in\Re} \int_{\St}
(x_j-\alpha)^2 \frac{1}{2}f_{t_i}(0)\vol_M(x_1,\ldots,x_n)\d{x_1}\cdots\d{x_n}\rb .
\]
Because the $x_j$ are smooth and orthonormal at $w$, we can choose
$\delta$ small enough so that the ratio of the Euclidean volume form
determined by the $x_j$ to the Riemannian volume form (which, we recall, includes the factor $H_0(x,\cdot)H_0(\cdot,y)$) is bounded from
below by $1/\sqrt{D}$ and from above by $\sqrt{D}$ on $B_\delta(w)$, for some positive constant $D$.  Assuming
that this is the case, we have that
\[
V(t_i) \geq \max_{j\in\{1,\ldots,n\}}\frac{1}{D}f_{t_i}(0)\lb \min_{\alpha\in\Re} \int_{\St}
(x_j-\alpha)^2 \d{x_1}\cdots\d{x_n}\rb
\]
and that the measure of $\St$ with respect to the Euclidean volume induced
by the $x_j$ (which we're now integrating against) is at least
$\pi_\delta(t_i\log2)/2$.

Now suppose we take the Steiner symmetrization (see
\cite{StroockIntro} for the definition and basic properties of Steiner
symmetrization) of $\St$ with respect to $x_1$.  Then
it's clear that the above integral with respect to $x_1$ is now
minimized by taking $\alpha=0$ and that it can only
have decreased after the symmetrization, while the integral with respect to
any other coordinate remains unchanged.  Hence, if we symmetrize with
respect to all of the coordinates and call the resulting set
$\Sym(\St)$, we see that
\[
V(t_i) \geq \max_{j\in\{1,\ldots,n\}} \frac{1}{D}f_{t_i}(0) \int_{\Sym(\St)} x_{j}^{2}
\d{x_1}\cdots\d{x_n} .
\]
Summing over $j$ and writing $r=x_1^2+\cdots+x_n^2$, we have that
\[
V(t_i) \geq \frac{1}{Dn} f_{t_i}(0)
\int_{\Sym(\St)} r^{2} \d{x_1}\cdots\d{x_n} .
\]
We recall that Steiner symmetrization preserves the Lebesgue measure
of sets, and thus the above integral is being taken over a set of
measure at least $\pi_{w,\delta}(t_i\log2)/2$.  On the other hand,
it's clear that for a given measure for $\Sym(\St)$,
this integral is minimized when $\Sym(\St)$ is the
ball of the appropriate measure centered at the origin.
Computing the integral of $r^2$ over such a ball, we conclude that
\be
\label{Eqn:Vt}
V(t_i) \geq \frac{\Area(\mS^{n-1})}{Dn^2(n+2)(2\omega_n)^{\frac{n+2}{n}}}
\lp\pi_{w,\delta}(t_i\log2)\rp^{\frac{n+2}{n}} f_{t_i}(0) .
\ee
Note that Equation~\eqref{Eqn:Vt} is valid (possibly with different constants) for any $w$ with associated
times $t_i$ (meaning that $\mu_{t_i}$ converges to a point mass at $w$) and smooth coordinates
$x_j$ which are orthonormal at $w$, once we choose small
enough $\delta$ and small enough $t_i$. 

Our next task is to compare the asymptotics of $\pi_{w,\delta}(t_i\log2)$ and
$f_{t_i}(0)=1/Z(t_i)$ as $t_i\da 0$.  Let $\pi(s)$ be the analogue of $\pi_{w,\delta}(s)$ over all of $\GE$.  Then
for any large enough $\alpha>0$, we have (up to exponentially small error, which we ignore)
\[
Z(t)= \pi(\alpha)e^{-\alpha/t} + \frac{1}{t}\int_0^\alpha
\pi(s)e^{-s/t}\d{s} .
\]
Choose any $\beta>0$.  Then for small enough $t$,
\[
\begin{split}
\frac{1}{t}\int_0^\alpha \pi(s)e^{-s/t}\d{s} &=
\frac{1}{t}\int_0^{-\beta t\log t}
\pi(s)e^{-s/t}\d{s} + \frac{1}{t}\int_{-\beta t\log t}^\alpha
\pi(s)e^{-s/t}\d{s} \\
& \leq \pi(-\beta t\log t)\lp 1-t^{\beta}\rp +\pi(\alpha)\lp
t^\beta-e^{-\alpha/t} \rp .
\end{split}
\]
It follows that, for small enough $t$, we have
\be
\label{Eqn:Zt}
Z(t) \leq \frac{1}{2}\pi(-\beta t\log t)+\pi(\alpha)t^\beta .
\ee

Next, we claim that $\pi(s)$ is bounded from above and below by positive powers of
$s$ (times a constant) for small enough $s$.  To see that, first note that $\rho$ is bounded from above by some non-negative function with a single zero and a non-degenerate Hessian at that zero, because the Hessian of $\rho$ over all of $\GE$ is controlled from above by compactness.  The volume of sublevel sets of this comparison function is easily seen to be asymptotic to a positive power of $s$, and thus $\pi(s)$ is bounded from below by this power of $s$.  Next, we observe that the Hessian of $\rho$ is always non-degenerate in the radial direction (in polar coordinates around $x$) and that the second derivative in this direction is bounded from below by a positive constant, again by compactness.  Hence $\rho$ is bounded from below by some quadratic function of the radius alone.  This comparison function also has the volume of its sublevel sets asymptotic to a positive power of $s$, and this provides an upper bound for $\pi(s)$.

Now let $a$ be the
infimum of all positive reals $b$ such that, for some constant $c$ and
small enough $t$, $ct^b\leq \pi(t)$.
Choose any small $\ol{\epsilon}>0$.  Then there must be some sequence
of times $\tau_i\da 0$ such that
\[
\frac{\pi(\tau_i)}{\tau_i^{a-\ol{\epsilon}}}\da 0 .
\]
Let the $\beta$ from equation \ref{Eqn:Zt} be given by $\beta= a-\ol{\epsilon}$, and let
$\ol{\tau}_i$ be a sequence such that $\tau_i=-\beta \ol{\tau}_i
\log\ol{\tau}_i$ for all $i$.  Then, applying equation \eqref{Eqn:Zt},
we have that
\[
\begin{split}
Z(\ol{\tau}_i) & \leq \pi(\tau_i) +\pi(\alpha)\ol{\tau}_i^\beta \\
& \leq \tau_i^{a-\ol{\epsilon}} +\pi(\alpha)\ol{\tau}_i^{a-\ol{\epsilon}} ,
\end{split}
\]
which implies that, for some constant $C$, we have
\be
\label{Eqn:ZTau}
\frac{1}{Z(\ol{\tau}_i)} \geq \frac{1}{C \ol{\tau}_i^{a-\ol{\epsilon}}
|\log\ol{\tau}_i|^{a-\ol{\epsilon}}}.
\ee
By compactness, we can assume, after passing to a subsequence, that
$\mu_{\ol{\tau}_i}$ converges to some limit measure, which must be a
point mass.  Thus, we can take $t_i = \ol{\tau}_i$.

In the remainder of the proof, we will assume that $C$ is some appropriate constant,
the exact value of which may vary from appearance to appearance.
Using that $t_i=\ol{\tau}_i$ and applying Equations~\eqref{Eqn:Vt}
and~\eqref{Eqn:ZTau}, we have that
\[
V(t_i) \geq C
\frac{\pi_{w,\delta}(t_i\log2)^{\frac{n+2}{n}}}{t_i^{a-\ol{\epsilon}}|\log
t_i|^{a-\ol{\epsilon}}} .
\]
By the same argument as before, we know
that $\pi(t\log2)$ is the measure of the subset of $\GE$ where the
density of $\mu_t$ is at least half of its maximum value.  Then,
because $\mu_{t_i}$ converges to a point mass at $w$, we must have that
$2\pi_{w,\delta}(t_i\log2)\geq \pi(t_i\log2)$ for small enough $t_i$.
Thus $\pi_{w,\delta}(t_i\log2) \geq Ct_i^{a+\ol{\epsilon}}$ (for some
potentially different constant $C$), and we conclude that  
\[
V(t_i)\geq C \frac{t_i^{2\ol{\epsilon}}}{|\log t_i|^{a-\ol{\epsilon}}}
\pi_{w,\delta}(t_i\log2)^{2/n} .
\]

Finally, it remains to estimate the asymptotics of $\pi_{w,\delta}(t)$.  Because we're assuming that $\mu_{t_i}$ converges to a point mass at $w$, Lemma \ref{Lem:Degenerate} implies that the Hessian of $\rho$ must be degenerate at $w$.
Since the Hessian of $\rho$ is degenerate at $w$, its Newton diagram
must be dominated by the diagram of $x_1^4+\sum_{i=2}^{n}x_i^2$,
possibly after relabeling the coordinates.  Hence our discussion of Laplace asymptotics (in particular, Theorem \ref{THM:Laplace}) implies that $\pi_{w,\delta}(t)\geq C t^{\frac{n}{2}-D_n}$ where $D_n$ is some positive constant depending only on the dimension
$n$ of $M$.  If we choose $\ol{\epsilon} \leq D_n/2n$, then we see that
\[
V(t_i) \geq C \frac{t_i^{1-D_n/n}}{|\log t_i|^{a-\ol{\epsilon}}} .
\]
This proves that there is some point $w$ and an associated sequence of
times $t_i$ such that, for any smooth coordinates $x_j$ orthonormal
at $w$, at least one of the $x_j$ has variance that goes to zero
slower than $t_i$.

Because $y\not\in\Cut(w)$, we can choose vectors $A_1,\ldots,A_n$
in $T_yM$ such that the functions $\nabla_{A_j}
E(z,y)-\nabla_{A_j}E(w,y)$ are smooth coordinates around $w$ which are
orthonormal at $w$.  Thus the above argument applies, and we conclude
that for some $A_j$, the variance of
$\nabla_{A_j}E(z,y)-\nabla_{A_j}E(w,y)$ goes to zero slower than
$t_i$.  Normalize this $A_j$ to have length one, and call the result
$A$.  Then $\nabla_{A}E(z,y)$ differs from
$\nabla_{A_j}E(z,y)-\nabla_{A_j}E(w,y)$ only by an affine
transformation, and thus its variance differs only by multiplication
by some (positive) constant.  We conclude that $\nabla_{A}E(z,y)$ has
variance that goes to zero slower than $t_i$, and this completes the
proof of the theorem for smooth manifolds.

In the case when $M$ is real-analytic, the proof follows the same lines, but is somewhat simpler.  First of all, we know that the limit measure $\mu$ is unique.  If it is not a point mass the result follows just as before, only now we don't have to pass to a subsequence of time and thus our result holds for the limit, not just the limit supremum.

Now suppose that $\mu_t$ converges to  point mass at some point $w$.  Then we have, using Equation~\eqref{Eqn:Vt} and noting that we don't need to pass to a subsequence,
\[
V(t) \geq \frac{\Area(\mS^{n-1})}{4n^2(n+2)(2\omega_n)^{\frac{n+2}{n}}}
\lp\pi_{w,\delta}(t\log2)\rp^{\frac{n+2}{n}} f_{t}(0) .
\]
Again, the point is to compare the asymptotics of $\pi_\delta(t\log2)$ and
$f_{t}(0)=1/Z(t)$.  We observe that $1/Z$ is essentially determined by the Laplace transform of $\pi_\delta$ (see \cite{Vasiliev} for a discussion of this fact).  So we are comparing a function with its Laplace transform.  In the proof of Theorem~\ref{THM:Unique}, we showed that $Z(t)\sim c t^{\alpha} |\log t|^m$ for some $c>0$, some non-negative rational number $\alpha$, and some non-negative integer $m$.  Direct computation shows that any function of this form is asymptotically equivalent (up to multiplication by some positive real number $a$) to the reciprocal of its Laplace transform.  (That this is true for real-analytic functions but not necessarily for smooth function is one reason for the increased difficulty in that case.)  Thus we have that
\[
V(t) \geq \frac{a \Area(\mS^{n-1})}{4n^2(n+2)(2\omega_n)^{\frac{n+2}{n}}}
\lp\pi_{w,\delta}(t\log2)\rp^{\frac{2}{n}} .
\]
Again, we know that $\pi_{w,\delta}(t)\geq C t^{\frac{n}{2}-D_n}$.  From here, the theorem follows just as above.  $\Box$

\subsection{Lower order leading terms}

Theorem~\ref{THM:Char} shows that the Hessian of $E_t(x,y)$ always blows up on the cut locus.  However, we've already seen that, if the limiting measure is a point mass, there will be no $1/t$ term in the expansion.  In other words, the variance will go to zero, but it will to do more slowly than $t$.  In this case, it is more difficult to determine the leading term of the expansion of $\Hess E_t(x,y)$, since it involves further terms in the Laplace asymptotics.  Because of this, we can only discuss the simplest case.

We consider the case where $\Gamma$ consists of a single point, $z_1$, such that $g$ can be diagonalized in a neighborhood of $z_1$.  Recall this means that there exists coordinates $u_j$ around $z_1$ such that
\[
g(u_1,\ldots,u_n)=\sum_{j=1}^{n} u_j^{2k_j}
\]
for some positive integers $k_1\leq \cdots\leq k_n$.  Since we assume $y\in\Cut(x)$, this geodesic must be conjugate and thus $k_n\geq 2$.  Further, since the Hessian of $h_{x,y}$ is non-degenerate in the radial direction, $k_1=1$.  The advantage to this situation is that we have a full asymptotic expansion from Equation \eqref{Eqn:KEExpansionFull}, rather than just the leading term (at least that's the advantage of assuming diagonalizability; assuming that $\Gamma$ contains just a single point makes the computation tractable).  Let $l$ be the smallest index such that $k_l=k_n$.  (Note that we allow $l=n$, but we must have $l\geq 2$.)  Thus $u_l,\ldots,u_n$ correspond to the directions of maximal conjugacy.  Keeping the first two terms of the expansion in Equation \eqref{Eqn:KEExpansionFull}, we see that
\[
\begin{split}
\int_{B_{\epsilon}(z_1)} & \phi(u_1,\ldots,u_n) \exp\lb -\frac{\sum_{j=1}^n u_j^{2k_j}}{t}\rb \d{u_1}\cdots\d{u_n} \sim \\
 &\frac{\Gamma(1/2k_1)\cdots \Gamma(1/2k_n)}{k_1\cdots k_n} t^{1/2k_1+\cdots+1/2k_n} \phi(0) \\
&+ \sum_{j=l}^n \frac{\Gamma(1/2k_1)\cdots \Gamma(3/2k_j) \cdots \Gamma(1/2k_n)}{k_1\cdots (2k_j)\cdots k_n} t^{1/2k_1+\cdots 3/2k_j\cdots+1/2k_n} \frac{\partial^2}{\partial u_j^2} \phi(0) .
\end{split}
\]

We will be interested in the cases $\phi=(\nabla_A E(z,y))^{*} H_0(x,z)H_0(z,y) \vol_u(z)$, where $*$ denotes 1 or 2 (corresponding to the expectation of the square and the square of the expectation which appear in the variance of $\nabla_A E (z,y)$), and $\phi= H_0(x,z)H_0(z,y) \vol_u(z)$ (which gives the expansion of $Z(t)$).  Dividing the relevant expansions and keeping the first two terms shows that (here we use that $k_j=k_n$ for all $j=l,\ldots,n$)
\[
\begin{split}
\E^{\mu_t}\lb (\nabla_A E(z,y))^{*} \rb \sim& (\nabla_A E(0,y))^{*} 
+\frac{\Gamma(3/2k_n)}{2\Gamma(1/2k_n)}t^{1/k_n} \sum_{j=l}^n \lb \frac{\partial^2}{\partial u_j^2}  (\nabla_A E(0,y))^{*} \right. \\
&\left. +2\frac{\frac{\partial}{\partial u_j}  (\nabla_A E(0,y))^{*} \frac{\partial}{\partial u_j}\lp  H_0(x,0)H_0(0,y) \vol_u(0)\rp}{ H_0(x,0)H_0(0,y) \vol_u(0)}\rb .
\end{split}
\]
If we now compute the variance, most of this cancels, and we see that
\[
\Var^{\mu_t}\lp \nabla_A E(z,y)\rp \sim \frac{\Gamma(3/2k_n)}{\Gamma(1/2k_n)}t^{1/k_n} \sum_{j=l}^{n} \lp \frac{\partial}{\partial u_j}  \nabla_A E(0,y)\rp^2 .
\]
Finally, we compute that
\[
\Hess_{A,A} E_t(x,y) \sim -\frac{4}{t^{1-1/k_n}} \frac{\Gamma(3/2k_n)}{\Gamma(1/2k_n)} \sum_{j=l}^{n} \lp \frac{\partial}{\partial u_j}  \nabla_A E(0,y)\rp^2 .
\]

There are several things to observe regarding this formula.  First of all, it shows that it's easy to produce situations in which $\Hess_{A,A}E_t(x,y)$ blows up at a rate intermediate between $1/t$ and $1$, and that every rational of the form $-(m-1)/m$
for a positive integer $m$ can be achieved as the order of the leading term.  This gives concrete intuition to the results of Theorem \ref{THM:Char}.  Second, in the case of a single minimal geodesic, we see that knowing the order of the leading term tells us the maximum order of degeneracy of the Hessian of $h_{x,y}$ and thus also ``how conjugate'' this minimal geodesic is, at least in the ``most conjugate'' directions.  Finally, it's relatively easy to see that the coefficient of this leading term also has geometric significance.  In particular, as a function of $A\in T_yM$, it is a symmetric, non-positive definite quadratic form such that the dimension of its kernel is equal to $n-l$.  To see this, let $\gamma(z)$ be the minimal geodesic from $z\in B_{\epsilon}(z_1)$ to $y$ and let $v_j\in T_yM$ be the derivative of the unit tangent to $\gamma(z_1)$ at $y$ with respect to $\partial/\partial u_j$.  Then $\partial/\partial u_j \nabla_A \dist(z_1,y)$ is non-zero if and only if $\ip{v_j}{A}$ is non-zero.   Because $h_{x,y}$ is non-degenerate in the radial direction, we can assume that the $\partial/\partial u_j$ (for $j=l,\cdots,n$) are perpendicular to $\gamma(z_1)$, and thus we see that
\[
\frac{\partial}{\partial u_j}  \nabla_A E(0,y) = \dist(z_1,y) \frac{\partial}{\partial u_j}  \nabla_A \dist(z_1,y) \quad\text{for $j=l,\ldots,n$} .
\]
Observing that, for $j=l,\cdots,n$, the $v_j$ are linearly independent, the desired result follows.  Thus, the leading term not only gives us the maximum order of degeneracy of the Hessian of $h_{x,y}$, but also tells us the dimension of the subspace on which this maximum degeneracy is achieved.

We find it interesting that, at least in this case, a term in the expansion further down than $1/t$ has such a nice geometric interpretation.  Unfortunately, these terms are hard to compute in general (the present case of a single minimal geodesic is the easiest case, yet even here completely working out the above calculations is rather laborious), and we don't know anything more about them.

\section{Mollification of Energy}

In the last section, we studied the asymptotics relative to fixed points $x$ and $y$.  Now, we turn our attention to considering how Corollary~\ref{TheCor} can be used to study the distributional Hessian of $E(x,y)$, where we think of $x$ as a fixed base point and thus of $E(x,y)$ as a function of $y$.  In particular, Varadhan's result (see Equation~\eqref{Eqn:Varadhan}) implies that $E_t(x,y)$ is a smooth mollifier of $E(x,y)$ as $t \da 0$.  Hence, computing the distributional limit of $\Hess_{A,A} E_t(x,y)$ as $t\da 0$ gives $\Hess_{A,A}E(x,y)$ as a distribution.

We note that the approach below is not the only way to study the distributional Hessian of  $E(x,y)$.  For example, many of these results are consequences of the fact that the cut locus is rectifiable with respect to $(n-1)$-dimensional Hausdorff measure, as proven by 
Mennucci~\cite{Mennucci} using viscosity solutions of Hamilton-Jacobi equations and geometric measure theory.

\subsection{The results}

Away from $\Cut(x)$, the distribution $\Hess_{A,A}E(x,y)$ is just a smooth function, and Equation~\eqref{Eqn:SAndM} tells us that $\Hess_{A,A}E_t(x,y)$ converges to this limit uniformly on compact subsets of $M\setminus\Cut(x)$.  This means that the singular part of $\Hess_{A,A}E(x,y)$, which we denote $\sing(\Hess_{A,A}E(x,y))$, is supported on $\Cut(x)$.  Considering Theorem~\ref{THM:Hess}, we see that any contribution to the singular part must come from the variance term.  In particular, for any smooth function $\phi$ we have
\[
\ip{\phi}{\sing(\Hess_{A,A}E(x,\cdot))} = -\lim_{\epsilon\da 0}
\lim_{t\da 0} \int_{B_{\epsilon}\lp\Cut(x)\rp} \phi(y) \frac{4}{t}
\Var^{\mu_{t,y}}(\nabla_{A}E(\cdot,y)) \d{y}
\]
where $\mu_{t,y}$ is the measure $\mu_{t}$ from above corresponding
to the point $y$, the variance of $\nabla_{A}E(\cdot,y)$ is taken with respect to the first variable, and we use the notation $\ip{\phi}{D}$ to denote the action of the distribution $D$ on the smooth function $\phi$.

Let $(r,\theta)\in\Re_{+}\times\mS^{n-1}$ be (normal) polar coordinates
around $x$, and let $d(\theta)$ be the distance to the cut locus along the geodesic
corresponding to $\theta$.  Let $U_x=\{(r,\theta):\theta\in \mS^{n-1},
r\in(d(\theta)-\epsilon,\d(\theta))\}$.  Then the exponential map gives a
diffeomorphism from $U_x$ to $M\setminus\Cut(x)$, and $\partial{U_x}=(d(\theta),\theta)$ is the tangential cut locus, that is, the preimage of $\Cut(x)$ under the exponential map (or more accurately, the connected component of the preimage closest to the origin).  This gives a natural identification of $\mS^{n-1}$ with the set of minimal geodesics from $x$ to $\Cut(x)$ and with $\partial U_x$.  We will frequently assume this identification, for example, when stating that some
$\theta\in\mS^{n-1}$ corresponds to a conjugate geodesic.  Because $\Cut(x)$ has measure zero, we can write the above integral in polar coordinates on $U_x$ as
\begin{multline}
\label{Eqn:Sing}
\ip{\phi}{\sing(\Hess_{A,A}E(x,y))} = \\
- \lim_{\epsilon\da 0}\lim_{t\da 0}
\int_{\mS^{n-1}} \left[ \int_{d(\theta)-\epsilon}^{d(\theta)} \phi(r,\theta)
\frac{4}{t} \Var^{\mu_{t,(r,\theta)}}(\nabla_{A}E(z,(r,\theta)))
\vol(r,\theta) \d{r} \right] \d{\theta}.
\end{multline}
We are now in a position to state the following theorem.
\begin{THM}
\label{THM:rho}
Let $M$ be a smooth, compact Riemannian manifold and let $x$ be any
point in $M$.  Let $A$ be any smooth vector field on $M$.  Choose
(normal) polar coordinates on $T_{x}M$ and define $U_x$ as above.
Then the right-hand side of
equation~\eqref{Eqn:Sing} defines a negative measure on $\partial U_x$, which
is absolutely continuous with respect to the measure $\d{\theta}$ on
$\partial U_x$ obtained by identifying it with $\mS^{n-1}$ via polar
coordinates.  Denote the corresponding Radon-Nikodym derivative by
$\rho(\theta)$; then $\rho(\theta)$ is bounded.  Thought of as a distribution on $M$,
$\Hess_{A,A}E(x,y))$ has as its singular part a negative measure $\nu_{x,A}$
supported on $\Cut(x)$, and further, $\nu_{x,A}$ is given by the
pushforward of $\rho(\theta)\d{\theta}$ under the exponential map.
\end{THM}

While Theorem~\ref{THM:rho} shows that the singular part of $\Hess_{A,A}E(x,y))$ has a relatively nice structure, it says little about the relationship between $\rho$ and the geodesic geometry of $M$.  In order to describe the relationship, we will need a bit more notation.  Let $C\subset \mS^{n-1}$ be the set of all $\theta$
which correspond to conjugate geodesics.  Next, say that the geodesics
corresponding to $\theta$ and $\tilde{\theta}$
are \emph{associated} if they lead to the same point in
$\Cut(x)$ (that is, if $d(\theta)=d(\tilde{\theta})$ and $(d(\theta),\theta)$ and
$(d(\tilde{\theta}),\tilde{\theta})$ are mapped to the same point under
$\exp_x$).  Let $P\subset \mS^{n-1}$ be the set of
$\theta\in\mS^{n-1}\setminus C$ to which there is associated
precisely one other $\tilde{\theta}\in \mS^{n-1}$ and such that
$\tilde{\theta}\not\in C$.  Finally, let $R=
\mS^{n-1}\setminus(C\cup P)$ (so $R$ consists of non-conjugate
$\theta$ which are associated to more than one other geodesic or which
are associated to a conjugate geodesic).  The three sets $C$, $P$ and
$R$ are disjoint and partition $\mS^{n-1}$.
\begin{THM}
\label{THM:CPR}
Let the hypotheses be as in Theorem~\ref{THM:rho}.
If $\theta\in C$, then $\rho(\theta)=0$.
Also, $R$ has
measure zero as a subset of $\mS^{n-1}$ with respect to $\d{\theta}$,
and $\rho$ is continuous except possibly at points of $R$.  Finally, there is an explicit expression for $\rho$ on $P$ (see Equation~\eqref{Eqn:RhoEquals} below).
\end{THM}

In order to give the explicit expression for $\rho$ on $P$,
we will need still more notation.  Let
$\theta$ be in $P$, let $\tilde{\theta}$ be the (one) associated
geodesic, and let $y$ be their common endpoint.  Also, let $z$ be the midpoint of the geodesic corresponding to
$\theta$, and let $B$ be the (non-degenerate) Hessian of $h_{x,y}$ at $z$.  Let $\tilde{z}$ and $\tilde{B}$ be the corresponding objects
associated to $\tilde{\theta}$.  Next, let $A_y\in T_yM$ be the value
of the vector field $A$ at $y$.  Then let $\psi$ be the angle
between the geodesic given by $\theta$ and $A_y$, $\tilde{\psi}$ be
the angle between the geodesic corresponding to $\tilde{\theta}$ and
$A_y$, and $\phi$ the angle between the geodesics $\theta$ and
$\tilde{\theta}$.  Then for any $\theta\in P$,
\begin{multline}\label{Eqn:RhoEquals}
\rho(\theta)= -\dist(x,y)|A_y|^2 \lp\cos\psi-\cos\tilde{\psi}\rp^2
  \vol(d(\theta),\theta) \\
  \times\lb (1-\cos\phi) \lp1 +\frac{H_0(x,z)H_0(y,z)\sqrt{\det
  \tilde{B}}}{H_0(x,\tilde{z})H_0(y,\tilde{z}) \sqrt{\det B}}\rp \rb^{-1}.
\end{multline}
Note that the volume element, all of the functions $H_0$
appearing above, and both $B$ and $\tilde{B}$ can be computed from the
Jacobi fields along the geodesics given by $\theta$ and $\tilde{\theta}$.

Theorem~\ref{THM:CPR} tells us that the only contributions to the
singular part come from points in $P$, which on $M$ are places where
locally the cut locus looks like a smooth
hypersurface and the singular part of $\Hess E(x,y)$ is just
given by the jump discontinuity of $\nabla E(x,y)$ across this
hypersurface.  While the cut locus itself can be quite complicated (for
example, it may not be triangulable, as shown by Gluck and Singer~\cite{GS}), the singular part of $\Hess E(x,y)$ is supported only at those points with the nicest local structure.

We should point out, however, that even though there
may not by any singular part of $\Hess E(x,y)$ in a neighborhood of
a conjugate point, the Hessian will not be smooth at a conjugate point.  To be precise, suppose that $\gamma:[0,\dist(x,y)]\mapsto M$ is a minimal geodesic (with unit speed parametrization) from $x$ to $y$.  Then if we consider $\Hess E(x,\gamma(s))$ as $s\nearrow \dist(x,y)$, we have that the Hessian, as an operator on vectors fields $A$, will blow up if and only if $\gamma$ is a conjugate geodesic.  Further, whether or not this blow up occurs for a given $A$ depends on how $A$ relates to the directions in which $\gamma$ is conjugate (we give the precise formula in the next section).

This leads us to the following picture.  The distributional Hessian of $E(x,y)$ is composed of two pieces, an $L^1$ function, which is just the regular Hessian on $M\setminus \Cut(x)$, and a singular part, which is the measure supported on $\Cut(x)$ as described above.  At a point $y$ which is the image on $M$ of a point in $P$, the $L^1$ part stays bounded as we approach $y$ along either minimal geodesic, but the singular measure is supported at $y$, as discussed above.  On the other hand, if $y$ is a conjugate point, then the singular measure may not be supported at $y$, but the $L^1$ part will blow up as we approach $y$ along any conjugate geodesic.  Because the set $R$ has measure zero in the decomposition of the tangential cut locus, these two cases completely describe the non-smooth behavior of the distributional Hessian of $E(x,y)$.

\subsection{The proofs}

The proof of Theorems~\ref{THM:rho} and~\ref{THM:CPR} will require a little preparation.  The first thing we need to do is to justify exchanging the integration with the limits in Equation~\eqref{Eqn:Sing}.  Equation~\eqref{Eqn:Sing} immediately implies that the singular part is non-positive, and thus is a non-positive measure (as opposed to a higher order distribution).  Because the singular part is a non-positive measure, we can estimate the Hessian in terms of the Laplacian.  In particular, choose any closed, connected set $\Omega\subset \mS^{n-1}$ the boundary of which has finite $(n-2)$-dimensional Hausdorff measure, and let $\Omega(\epsilon)$ be the set
\[
\left\{(r,\theta) : \theta\in\Omega\text{ and }r\in [d(\theta)-\epsilon,d(\theta)] \right\} .
\]
Then we have
\begin{multline*}
\left| \int_{\Omega(\epsilon)}  \phi(r,\theta)
\Hess_{A,A} E_t(x,y) \vol(r,\theta) \d{r} \d{\theta}\right|
 \leq C_1 |\phi|_{\infty}|A|^2 \\
\times\lc \left| \int_{\Omega(\epsilon)} \Lap E_t(x,y)
\vol(r,\theta) \d{r}  \d{\theta}  \right| +
\int_{\Omega(\epsilon)}  \vol(r,\theta) \d{r}  \d{\theta}  \rc 
\end{multline*}
for small enough $t$ and $\epsilon$ and some constant $C_1$.  Here $|f|_{\infty}$ denotes the $L^{\infty}$ norm of $f$ on the set $\Omega(\epsilon)$.  This last integral can be estimated as
\[
\int_{\Omega(\epsilon)}  \vol(r,\theta) \d{r}  \d{\theta} \leq
\epsilon  |\Omega|_{\mS^{n-1}} |\vol(r,\theta)|_{\infty} .
\]

Let $n(r,\theta)$ be the outward pointing unit normal to $\Omega(\epsilon)$.  We use integration by parts to write
\[
\int_{\Omega(\epsilon)} \Lap E_t(x,y) \vol(r,\theta) \d{r}\d{\theta} = 
- \int_{\partial \Omega(\epsilon)} \ip{\nabla E_t(x,y)}{n(r,\theta)} \d{\sH^{n-1}_{M}} ,
\]
where $\sH^{n-1}_{M}$ denotes $(n-1)$-dimensional Hausdorff measure relative to the Riemannian volume measure.  The right-hand side is a smooth function of $t$, so the only question is what happens to it as $t\da 0$.  On the ``walls'' of $\partial \Omega(\epsilon)$, that is, the set
\[
W(\Omega(\epsilon)) =\left\{(r,\theta) : \theta\in\partial\Omega\text{ and }r\in (d(\theta)-\epsilon,d(\theta)) \right\} ,
\]
the results of Stroock and Malliavin show that $\nabla E_t(x,y)$ converges (pointwise) to $\nabla E(x,y)=r\partial_r$.  Further, the Gauss lemma shows that, on this set, $n(r,\theta)$ is perpendicular to $\partial_r$ (both in the Riemannian metric and the Euclidean metric).  Thus, the integral over the ``walls'' goes to zero with $t$, which we will write
\[
 \int_{W(\Omega(\epsilon))} \ip{\nabla E_t(x,y)}{n(r,\theta)} \d{\sH^{n-1}_{M}} =o(1) .
\]
Next, we consider the integral over the ``top'' and ``bottom'' of $\partial \Omega(\epsilon)$, that is, the sets
\[ \begin{split}
T(\Omega(\epsilon)) =&\left\{(r,\theta) : \theta\in\partial\Omega\text{ and }r=d(\theta) \right\} \\
\text{and}\quad B(\Omega(\epsilon)) =&\left\{(r,\theta) : \theta\in\partial\Omega\text{ and }r=d(\theta)-\epsilon \right\} .
\end{split} \]
The results of Stroock and Turetsky show that the norm of $\nabla E_t(x,y)$ is bounded for all small $t$, and thus  $\ip{\nabla E_t(x,y)}{n(r,\theta)}$ is bounded on both $T(\Omega(\epsilon))$ and $B(\Omega(\epsilon))$.  Thus it remains only to control the $\sH^{n-1}_{M}$ measure of $T(\Omega(\epsilon))$ and $B(\Omega(\epsilon))$.

First note that the volume density $\vol(r,\theta)$ is bounded, and thus we can estimate the $\sH^{n-1}_{M}$ measure of a set from above by $|\vol(r,\theta)|_{\infty}$ times the $\sH^{n-1}$ measure of the set, where we use $\sH^{n-1}$ to indicate the $(n-1)$-dimensional Hausdorff measure relative to the Euclidean volume on the tangent space.  The key to estimating this measure is the result of Itoh and Tanaka~\cite{ItohTanaka} that the distance to the cut locus (that is, the function $d(\theta)$) is Lipschitz.  It follows that
\[
\sH^{n-1}\lp T(\Omega(\epsilon))\cup B(\Omega(\epsilon)) \rp \leq 2\lp 1+\Lip(d)\rp |\Omega|_{\mS^{n-1}} |d|_{\infty} .
\]
Combining the above results, we conclude that
\begin{multline}
\label{Eqn:BasicEstimate}
\left| \int_{\Omega} \left[ \int_{d(\theta)-\epsilon}^{d(\theta)} \phi(r,\theta)
\Hess_{A,A} E_t(x,y)
\vol(r,\theta) \d{r} \right] \d{\theta}\right|
 \leq \\
C_2 |\phi|_{\infty} |A|^2 |\vol(r,\theta)|_{\infty} |\Omega|_{\mS^{n-1}}  \lb \lp 1+\Lip(d)\rp |d|_{\infty}+o(1)  + \epsilon  \rb
\end{multline}
for small enough $t$ and $\epsilon$ and some constant $C_2$.

This justifies exchanging the limit as $t\da 0$ with the integral over $\theta$ in equation \eqref{Eqn:Sing}.  In addition, because the variance and volume density are always non-negative and $\phi$ is smooth, we can also exchange the limit as $\epsilon\da 0$ with the integral over $\theta$.  This shows that the measure given by the right-hand side of equation \eqref{Eqn:Sing} is absolutely continuous with respect to the measure $\d{\theta}$ and that its Radon-Nikodym derivative $\rho$ is given by
\be
\label{Eqn:rho}
\rho(\theta) = -\lim_{\epsilon\da 0}\lim_{t\da 0}
\int_{d(\theta)-\epsilon}^{d(\theta)}
\frac{4}{t} \Var^{\mu_{t,(r,\theta)}}(\nabla_{A}E(z,(r,\theta)))
\vol(r,\theta) \d{r}  .
\ee

We now turn our attention to the sets $C$, $R$, and $P$.  On the set $C$, we want to show that $\rho=0$.  To do so, we won't work directly with \eqref{Eqn:Sing}, but rather with the definition of the distributional Hessian.  In particular, choose any $\theta\in C$, and let $\Omega_{\delta}\subset\mS^{n-1}$ be a disk of radius $\delta$ around $\theta$.  Then starting from \eqref{Eqn:BasicEstimate} and letting $t$ and $\epsilon$ both go to zero, we see that
\[
- \int_{\Omega_{\delta}} \rho(\psi)\d{\psi} \leq C_2 |A|^2 |\vol(d(\psi),\psi)|_{\infty} |\Omega_{\delta}|_{\mS^{n-1}}  \lp 1+\Lip(d)\rp |d|_{\infty} .
\]
Everything on the right-hand side is bounded from above, so if we divide both sides by $|\Omega_{\delta}|_{\mS^{n-1}}$ and let $\delta$ go to zero, then Lebesgue's differentiation theorem tells us that $-\rho(\theta)$ is almost everywhere on $C$ less than or equal to a constant times the limit as $\delta\da 0$ of the $L^{\infty}$-norm of $\vol(d(\psi),\psi)$ over $\Omega_{\delta}$.  Because $\vol(d(\psi),\psi)$ is continuous and equals zero at $\theta$ (because the corresponding geodesic is conjugate), this shows that $\rho$ is equal to zero almost everywhere on $C$.  Because $\rho$ is a density, its value only matters up to almost everywhere equivalence, and thus we can take $\rho$ to be zero on all of $C$.

Next, we prove a lemma which shows that when computing $\rho$, we can ignore the set $R$.  
\begin{Lemma}
For any compact manifold $M$ and a basepoint $x\in M$, the corresponding set $R$, as defined above, has measure zero as a subset of $\mS^{n-1}$ with its standard volume measure.
\end{Lemma}
\emph{Proof.}  As usual, we identify $\mS^{n-1}$ with the set of minimal geodesics from $x$ to $\Cut(x)$.  Let $S_1$ be set set of non-conjugate $\theta$ which are associated to more than one other geodesics, all of which are non-conjugate.  Let $S_2$ be the set of non-conjugate $\theta$ which are associated to a conjugate geodesic (and possibly to other geodesics as well).  Then $R=S_1\cup S_2$.

We first consider $S_1$.  In particular, choose any $\theta\in S_1$ and let and let $\theta_1,\ldots,\theta_k$ for some $k\geq 2$ be the associated geodesics (there are necessarily only finitely many because otherwise they would have an accumulation point, forcing at least one to be conjugate).  Let $y$ be their common endpoint.  In terms of the exponential map, the fact that $\theta$ is not conjugate means that there is a neighborhood $U$ of $(d(\theta,\theta)$ in $T_xM$ which is diffeomorphic to a neighborhood $V$ of $y$ under the exponential map.  The same is true for each of the $\theta_i$ and we use $U_i$ for the corresponding subsets of $T_xM$.  By choosing these neighborhoods small enough, we can assume that the preimage of $V$ under the exponential map is precisely equal to $U\cup U_1 \cup\cdots \cup U_k$.  Now we can define a smooth function $f$ on $V$ to be the length of the corresponding element of the tangent space in $U$, and we can similarly define the $f_i$.  We will call such functions local distance functions.  It follows that, for any point in $V$, the Riemannian distance is given by $\min\{f,f_1,\ldots,f_k\}$ and further, that the number of minimal geodesics to that point is given by the number of these local distance functions which achieve this minimum.  At $y$, we know that all of these local distance functions achieve their common minimum, that they all have gradient of length one, and that these gradient vectors are all distinct.  Given this, it is straight-forward to see that the set of points in $V$ where at least three of the these local distance functions achieve the common minimum is given by a finite union of smooth submanifolds of dimension no more than $n-2$.  Since $V$ and $U$ are diffeomorphic under the exponential map, it follows that the set of points in $U$ which correspond to points with more than two minimal geodesics has $(n-1)$-dimensional Hausdorff measure equal to zero.  Projecting this set onto the $\theta$ coordinate gives the intersection of $S_1$ with a neighborhood of $\theta$, and we conclude that this set has measure zero.  This shows that every $\theta\in S_1$ has a neighborhood in which $S_1$ has measure zero, and it follows that all of $S_1$ has measure zero.

We now consider $S_2$.  This will require some facts about the set of points in $M$ which are conjugate to $x$ along a minimal geodesic; call this set the conjugate-cut locus.  Observe that the set of points in $M$ corresponding to $S_2$ is contained in the conjugate-cut locus.  Let $C_1$ be the set of all vectors in $T_xM$ such that the kernel of $d\exp_x$ has dimension $1$.  In~\cite{Warner1}, Warner showed that $C_1$ is a smooth $(n-1)$-dimensional submanifold of $T_xM$.  Let $T$ be the set of points in $C_1$ such that the kernel of $d\exp_x$ is contained in the tangent space to $C_1$ at that point.  Let $H$ be the set of all vectors in $T_xM$ such that the kernel of $d\exp_x$ has dimension at least two.  In the proof of Lemma~1.1 of~\cite{Warner2}, Warner showed that the image under the exponential map of $T\cup H$ has $(n-1)$-dimensional Hausdorff measure equal to zero.  In the proof of Proposition~3.2 of~\cite{Hebda}, Hebda proved that a point of $C_1\setminus T$ cannot correspond to a minimal geodesic.  This means that the conjugate-cut locus is precisely the image of $T\cup H$ under $\exp_x$ and thus has $(n-1)$-dimensional Hausdorff measure equal to zero.  Now choose any $\theta\in S_2$ and let $y$ be the corresponding point in $M$.  As before, we can choose neighborhoods $U$ and $V$ of $(d(\theta),\theta)$ and $y$ such that the exponential map gives a diffeomorphism between them.  In particular, this means that there is a neighborhood $W\subset \mS^{n-1}$ of $\theta$ such that no two points in $W$ correspond to the same point in $M$.  Suppose $S_2\cap W$ has positive measure.  Then the corresponding set of points in $U$ would have positive $(n-1)$-dimensional Hausdorff measure.  But these points are contained in the conjugate-cut locus and so this contradicts the above.  Hence we conclude that $S_2\cap W$ has measure zero.  Since this holds for any $\theta\in S_2$, it follows that all of $S_2$ has measure zero.  Given the decomposition $R=S_1\cup S_2$, this completes the proof of the lemma.  $\Box$

In order to complete the proofs of Theorems~\ref{THM:rho} and~\ref{THM:CPR}, it suffices to compute $\rho$ from equation~\eqref{Eqn:rho} on $P$.  This is fairly straight-forward, if somewhat tedious.  In this case, $\sO_{\epsilon}$ consists just of balls around $z$ and $\tz$.  We begin by estimating $\E^{\mu_{t,\delta}}\lb f\rb$ for any smooth test function $f$ and small $t$ and $\delta$.  First of all, let $y_{\delta}$ be the image under the exponential map of $(d(\theta)-\delta,\theta)$.  Then let $z_{\delta}$ be the midpoint of the geodesic from $x$ to $y_{\delta}$ in the direction $\theta$, and let $\tz_{\delta}$ be the midpoint of the (non-minimizing, if $\delta>0$) geodesic from $x$ to $y_{\delta}$ in the direction close to $\tilde{\theta}$.  Because the exponential map is a diffeomorphism near $(d(\theta),\tilde{\theta})$, it follows that $\tz_{\delta}$ is well-defined for small enough $\delta$ and depends smoothly on $\delta$.  We know that $h_{x,y_{\delta}}$ has non-degenerate Hessian at both $z_{\delta}$ and $\tz_{\delta}$, and we denote these Hessians by $B_{\delta}$ and $\tilde{B}_{\delta}$.  An easy computation shows that the volume form associated to coordinates which diagonalize $2h_{x,y_{\delta}}$ around $z_{\delta}$ is $1/\sqrt{\det B_{\delta}}$, and similarly for $\tz_{\delta}$; note that $B_{\delta}$ is the Hessian of $h_{x,y_{\delta}}$, without the factor of 2 (this is essentially a consequence of the $1/2$ which appears in the second order Taylor expansion).  Then Equation~\eqref{Eqn:KEExpansionFull} gives
\[\begin{split}
\int_{\sO_{\epsilon}} f(u) & H_0(x,u)H(u,y_{\delta})  \exp\lb -\frac{2h_{x,y_{\delta}}(u)}{t}\rb \d{u} =\\
& \Gamma\lp\frac{1}{2}\rp^n t^{n/2}\exp\lb -\frac{2h_{x,y_{\delta}}(z_{\delta})}{t}\rb\lb f(z_{\delta})\frac{H_0(x,z_{\delta})H_0(z_{\delta},y_{\delta})}{\sqrt{\det B_{\delta}}} +O_{\delta}(t) \rb \\
+& \Gamma\lp\frac{1}{2}\rp^n t^{n/2}\exp\lb -\frac{2h_{x,y_{\delta}}(\tz_{\delta})}{t}\rb\lb f(\tz_{\delta})\frac{H_0(x,\tz_{\delta})H_0(\tz_{\delta},y_{\delta})}{\sqrt{\det \tilde{B}_{\delta}}} +O_{\delta}(t) \rb . 
\end{split}\]
Here we use the notation $O_{\delta}(t)$ to indicate that the error term as $t\da 0$ depends on $\delta$, but does so uniformly for all sufficiently small $\delta$; in particular, the integral of $O_{\delta}(t)$ with respect to $\delta$ is $O(t)$.

Let
\[
F=\frac{H(x,\tz)H(y,\tz)\sqrt{\det B}}{H(x,z)H(y,z) \sqrt{\det \tilde{B}}} .
\]
Also note that $f(z_{\delta})=f(z)+O(\delta)$, $f(\tz_{\delta})=f(\tz)+O(\delta)$, and so on for $H_0$, $h$, and $B$.   In addition, we will need to know how $2h_{x,y_{\delta}}(z_{\delta})$ compares with $2h_{x,y_{\delta}}(\tz_{\delta})$.  Elementary trigonometry shows that this difference can be written as
\[
2h_{x,y_{\delta}}(z_{\delta})-2h_{x,y_{\delta}}(\tz_{\delta})= \dist(x,y)\delta(1-\cos\phi)+O(\delta^2) .
\]
Then since $Z(t)$ is obtained simply by taking $f(u)\equiv 1$ in the above, we have that
\begin{multline*}
\E^{\mu_{t,\delta}}\lb f\rb = \\
\frac{f(z)+O(\delta)+\lp f(\tz)F +O(\delta)\rp \exp\lb -\frac{1}{t} (\delta\dist(x,y)(1-\cos \phi) +O(\delta^2))\rb}{1+\lp F+O(\delta)\rp \exp\lb -\frac{1}{t} (\delta\dist(x,y)(1-\cos \phi) +O(\delta^2))\rb} +O_{\delta}(t) .
\end{multline*}
Using this, we can compute (after doing some algebra) that
\begin{multline*}
\Var^{\mu_{t,\delta}}\lb \nabla_A E(\cdot,y) \rb = \lp \nabla_A E(z,y)-\nabla_A E(\tz,y)\rp^2 \lp F+O(\delta)\rp \\
\times \frac{\exp\lb -\frac{1}{t} (\delta\dist(x,y)(1-\cos \phi) +O(\delta^2))\rb}{\left\{1+\lp F+O(\delta)\rp \exp\lb -\frac{1}{t} (\delta\dist(x,y)(1-\cos \phi) +O(\delta^2))\rb\right\}^2} +O_{\delta}(t).
\end{multline*}
Writing $\vol(d(\theta)-\delta,\theta)$ as $\vol(d(\theta),\theta)+O(\delta)$ and making the change of variables $\alpha = -\delta\dist(x,y)(1-\cos\phi)/t$, Equation \eqref{Eqn:rho} gives
\begin{multline*}
\rho(\theta)= \frac{4\lp \nabla_A E(z,y)-\nabla_A E(\tz,y)\rp^2}{\dist(x,y)(1-\cos\phi)} \vol(d(\theta),\theta) \times \\
 \lim_{\epsilon\da0}\lim_{t\da0} \int_{I(\epsilon,t)} \lb \frac{\lp F+O(t\alpha) \rp e^{\alpha+O((t\alpha)^2)}}{\lc 1+\lp F+O(t\alpha)\rp e^{\alpha+O((t\alpha)^2)}\rc^{2}} +O_{\delta}(t)\rb \lb 1+O(t\alpha) \rb \d{\alpha}
\end{multline*}
where $I(\epsilon,t)$ is the interval $[-\epsilon\dist(x,y)(1-\cos\phi)/t,0]$.

Taking both limits causes the region of integration to become $\lp -\infty,0 \rb$.  Further, it causes all of the $O(\cdot)$ terms to vanish.  To see this first note that $O_{\delta}(t)$ vanishes uniformly with $t$.  As for the $O(t\alpha)$ and $O((t\alpha)^2)$ terms, $O(t\alpha)$ is bounded on $I(\epsilon,t)$ for all $\epsilon$ and $t$ and goes to zero uniformly on any compact subinterval.  Using this, one can show that they don't contribute in the limit (to see this in detail, one can make the further change of variables $\beta= F e^{\alpha}$ and compute the integral).  Thus, taking the limits, the above integral becomes
\[
\int_{-\infty}^{0} \frac{F e^{\alpha}}{\lc 1+F e^{\alpha}\rc^{2}}\d{\alpha} = \frac{1}{1+F^{-1}}.
\]
Using this in the above expression for $\rho$, along with the fact that
\[
\lp \nabla_A E(z,y)-\nabla_A E(\tz,y)\rp^2 = \frac{1}{4}\dist(x,y)^2 \lab A_y \rab^{2}\lp \cos\psi-\cos\tilde{\psi} \rp^2 ,
\]
gives Equation~\eqref{Eqn:RhoEquals}.  This completes the proofs of Theorems~\ref{THM:rho} and~\ref{THM:CPR}.

Finally, we justify our earlier comments about the $L^1$ part of the distributional Hessian along conjugate geodesics.  We begin by applying Theorem~\ref{THM:Hess} at a point $y\not\in \Cut(x)$.  Because $y$ is not in the cut locus, there is a single, non-conjugate minimal geodesic between $x$ and $y$.  Because this geodesic is not conjugate, $\Hess 2h_{x,y}(z)$ is non-degenerate, and we choose coordinates $u_1,\ldots,u_n$ around $z$ such that $2h_{x,y} = u_1^2 +\cdots +u_n^2$.  We already know that, since we are not on the cut locus, the leading term in the Hessian will be the constant term, and in order to compute this term we will need the first two terms in the expansion of Theorem~\ref{THM:Hess}. We begin by computing, using the first two terms of the expansion in Equation~\eqref{Eqn:KEExpansionFull}, that
\begin{multline} \label{Eqn:BigExp}
\int_{B_\epsilon(z)} f(u) \exp \lb -\frac{2h_{x,y}(u)}{t}\rb k(t/2,x,u)k(t/2,y,u) \d{u} = \\
 \exp \lb -\frac{2h_{x,y}(z)}{t}\rb t^{n/2}  \Big\{ \Gamma(1/2)^n f(z) k(t/2,x,z)k(t/2,y,z)\vol_u(z)  \\
\left. + t \frac{\Gamma(1/2)^{n-1}\Gamma(3/2)}{2} \Lap^u \lp f(z) k(t/2,x,z)k(t/2,y,z)\vol_u(z) \rp +O(t^{2}) \rc .
\end{multline}
Here the symbol $\Lap^u$ means the operator $\sum_{i=1}^{n}\frac{\partial^2}{\partial u_i^2}$.  Recall that 
\[ \begin{split}
k(t,x,y) &=H_0(x,y)+tH_1(x,y)+O(t^2) \\
\text{and}\quad l(t,x,y,A) &= -\nabla_A E(x,y) +t\nabla_A G_1(x,y)+O(t^2) .
\end{split} \]

Given Equation~\eqref{Eqn:BigExp} and the expansions of $l$ and $k$ in terms of the $H_i$ and the $G_i$, expanding the right hand side of the equality  in Theorem~\ref{THM:Hess} becomes simply a lengthy exercise in manipulating Taylor series.  We won't reproduce the computation here and will instead merely state the result.  The coefficient of the $1/t$ term is zero, as we know it must be, and taking the limit as $t\da 0$ gives
\[
\Hess_{A,A} E(x,y) = 2\lb \Hess_{A,A} E(z,y) -\sum_{i=1}^{n} \lp\frac{\partial}{\partial u_i} \nabla_A E(z,y) \rp^2 \rb .
\]
If we consider what happens as $y$ approaches $\Cut(x)$ along a geodesic $\gamma$, we see that $\Hess E(x,y)$ will blow up, as an operator, if and only if $\gamma$ is conjugate.  This is because both $\Hess E(z,y)$ and $\nabla E(z,y)$ remain bounded, and thus the only way for a blow up to occur is if at least one of the $\partial_{u_i}$ has its (Riemannian) length blowing up.  This occurs precisely if the corresponding eigenvalue of $\Hess h_{x,y}(z)$ is going to zero, and thus precisely if $\gamma$ is conjugate in the direction corresponding to $\partial_{u_i}$.  It is also easy to see that the blow up must be in the negative direction, and that the relationship between a given vector $A$ and the $\partial_{u_i}$ determines whether or not the Hessian blows up for a given $A$.

\def\cprime{$'$} \def\cprime{$'$}
\providecommand{\bysame}{\leavevmode\hbox to3em{\hrulefill}\thinspace}
\providecommand{\MR}{\relax\ifhmode\unskip\space\fi MR }
% \MRhref is called by the amsart/book/proc definition of \MR.
\providecommand{\MRhref}[2]{%
  \href{http://www.ams.org/mathscinet-getitem?mr=#1}{#2}
}
\providecommand{\href}[2]{#2}


\begin{thebibliography}{10}

\bibitem{Arnold}
V.~I. Arnol{\cprime}d, S.~M. Guse\u{\i}~n Zade, and A.~N. Varchenko,
  \emph{Singularities of differentiable maps. {V}ol. {II}}, Monographs in
  Mathematics, vol.~83, Birkh\"auser Boston Inc., Boston, MA, 1988, Monodromy
  and asymptotics of integrals, Translated from the Russian by Hugh Porteous,
  Translation revised by the authors and James Montaldi.

\bibitem{BenArous}
G.~Ben~Arous, \emph{D\'eveloppement asymptotique du noyau de la chaleur
  hypoelliptique hors du cut-locus}, Ann. Sci. \'Ecole Norm. Sup. (4)
  \textbf{21} (1988), no.~3, 307--331.

\bibitem{BGV}
Nicole Berline, Ezra Getzler, and Mich{\`e}le Vergne, \emph{Heat kernels and
  {D}irac operators}, Grundlehren der Mathematischen Wissenschaften
  [Fundamental Principles of Mathematical Sciences], vol. 298, Springer-Verlag,
  Berlin, 1992.

\bibitem{BierstoneMilman}
Edward Bierstone and Pierre~D. Milman, \emph{Canonical desingularization in
  characteristic zero by blowing up the maximum strata of a local invariant},
  Invent. Math. \textbf{128} (1997), no.~2, 207--302.

\bibitem{Bishop}
Richard~L. Bishop, \emph{Decomposition of cut loci}, Proc. Amer. Math. Soc.
  \textbf{65} (1977), no.~1, 133--136.

\bibitem{Chavel}
Isaac Chavel, \emph{Eigenvalues in {R}iemannian geometry}, Pure and Applied
  Mathematics, vol. 115, Academic Press Inc., Orlando, FL, 1984, Including a
  chapter by Burton Randol, With an appendix by Jozef Dodziuk.

\bibitem{EAndK}
Ricardo Estrada and Ram~P. Kanwal, \emph{A distributional approach to
  asymptotics}, second ed., Birkh\"auser Advanced Texts: Basler Lehrb\"ucher.
  [Birkh\"auser Advanced Texts: Basel Textbooks], Birkh\"auser Boston Inc.,
  Boston, MA, 2002, Theory and applications.

\bibitem{GS}
Herman Gluck and David Singer, \emph{Deformations of geodesic fields}, Bull.
  Amer. Math. Soc. \textbf{82} (1976), no.~4, 571--574.

\bibitem{Hebda}
James~J. Hebda, \emph{The local homology of cut loci in {R}iemannian
  manifolds}, T\^ohoku Math. J. (2) \textbf{35} (1983), no.~1, 45--52.

\bibitem{Hsu}
Elton~P. Hsu, \emph{Stochastic analysis on manifolds}, Graduate Studies in
  Mathematics, vol.~38, American Mathematical Society, Providence, RI, 2002.

\bibitem{ItohTanaka}
Jin-ichi Itoh and Minoru Tanaka, \emph{The {L}ipschitz continuity of the
  distance function to the cut locus}, Trans. Amer. Math. Soc. \textbf{353}
  (2001), no.~1, 21--40.

\bibitem{SAndM}
Paul Malliavin and Daniel~W. Stroock, \emph{Short time behavior of the heat
  kernel and its logarithmic derivatives}, J. Differential Geom. \textbf{44}
  (1996), no.~3, 550--570.

\bibitem{Mennucci}
Andrea C.~G. Mennucci, \emph{Regularity and variationality of solutions to
  {H}amilton-{J}acobi equations. {I}. {R}egularity}, ESAIM Control Optim. Calc.
  Var. \textbf{10} (2004), no.~3, 426--451 (electronic).

\bibitem{MP}
S.~Minakshisundaram and {\AA}.~Pleijel, \emph{Some properties of the
  eigenfunctions of the {L}aplace-operator on {R}iemannian manifolds}, Canadian
  J. Math. \textbf{1} (1949), 242--256.

\bibitem{Molchanov}
S.~A. Mol{\v{c}}anov, \emph{Diffusion processes, and {R}iemannian geometry},
  Uspehi Mat. Nauk \textbf{30} (1975), no.~1(181), 3--59.

\bibitem{MeStroock}
Robert Neel and Daniel Stroock, \emph{Analysis of the cut locus via the heat
  kernel}, Surveys in differential geometry. Vol. IX, Surv. Differ. Geom., IX,
  Int. Press, Somerville, MA, 2004, pp.~337--349.

\bibitem{StroockIntro}
Daniel~W. Stroock, \emph{A concise introduction to the theory of integration},
  third ed., Birkh\"auser Boston Inc., Boston, MA, 1999.

\bibitem{SAndT2}
Daniel~W. Stroock and James Turetsky, \emph{Short time behavior of logarithmic
  derivatives of the heat kernel}, Asian J. Math. \textbf{1} (1997), no.~1,
  17--33.

\bibitem{SAndT}
\bysame, \emph{Upper bounds on derivatives of the logarithm of the heat
  kernel}, Comm. Anal. Geom. \textbf{6} (1998), no.~4, 669--685.

\bibitem{Sussmann}
H.~J. Sussmann, \emph{Real analytic desingularization and subanalytic sets: an
  elementary approach}, Trans. Amer. Math. Soc. \textbf{317} (1990), no.~2,
  417--461.

\bibitem{Varchenko}
A.~N. Var{\v{c}}enko, \emph{Newton polyhedra and estimates of oscillatory
  integrals}, Funkcional. Anal. i Prilo\v zen. \textbf{10} (1976), no.~3,
  13--38.

\bibitem{Vasiliev}
B.~A. Vasil{\cprime}ev, \emph{The asymptotic behavior of exponential integrals,
  the {N}ewton diagram and the classification of minima}, Funkcional. Anal. i
  Prilo\v zen. \textbf{11} (1977), no.~3, 1--11, 96.

\bibitem{Warner1}
Frank~W. Warner, \emph{The conjugate locus of a {R}iemannian manifold}, Amer.
  J. Math. \textbf{87} (1965), 575--604.

\bibitem{Warner2}
\bysame, \emph{Conjugate loci of constant order}, Ann. of Math. (2) \textbf{86}
  (1967), 192--212.

\end{thebibliography}
\end{document}